\documentclass[envcountsect,a4paper]{llncs}

\usepackage[latin1]{inputenc}   
\usepackage[dvips]{graphicx}
\usepackage{latexsym}
\usepackage{amssymb}
\usepackage{amsmath}

\newtheorem{lemma2}[theorem]{Lemma}
\newtheorem{cor}[theorem]{Corollary}
\newtheorem{proposition2}[theorem]{Proposition}

\title{Independence Complexes of Cylinders Constructed from Square and Hexagonal Grid Graphs}
\author{Johan Thapper\thanks{The author was supported by the 
    {\em Programme for Interdisciplinary Mathematics}, 
    Department of Mathematics, Link\"opings universitet}}
\institute{Department of Computer and Information Science\\
Link\"opings universitet\\
SE-581 83 Link\"oping, Sweden\\
\email{johth@ida.liu.se}}

\newcommand{\homeo}{\cong}
\newcommand{\homot}{\simeq}
\newcommand{\vv}{\mbox{\boldmath$v$}}

\usepackage{color}
\usepackage{epsfig}

\usepackage{multirow}

\begin{document}            
\maketitle

\begin{abstract}
\noindent
  In the paper [Fendley et al., J. Phys. A: Math. Gen., 38 (2005),
    pp. 315-322], Fendley, Schoutens and van Eerten studied
  the hard square model at negative activity.
  They found analytical and numerical evidence that the eigenvalues
  of the transfer matrix with periodic boundary were all
  roots of unity.
  They also conjectured that for an $m \times n$ square grid, with
  doubly periodic boundary, the partition function is equal to 1
  when $m$ and $n$ are relatively prime.
  These conjectures were proven in
  [Jonsson, Electronic J. Combin., 13(1) (2006), R67].
  There, it was also noted that the cylindrical case seemed to have
  interesting properties when the circumference of the cylinder
  is odd.
  In particular, when 3 is a divisor of both the circumference and
  the width of the cylinder minus 1, the partition function is -2.
  Otherwise, it is equal to 1.
  In this paper, we investigate the hard square and
  hard hexagon models at activity -1, with
  single periodic boundary, i.e, cylindrical identifications,
  using both topological and combinatorial techniques.
  We compute the homology groups of the associated independence complex
  for small sizes and suggest a matching which, we believe,
  with further analysis could help solve the conjecture.
  We also briefly review a technique recently described by
  Bousquet-Mélou, Linusson and Nevo, for determining some of the
  eigenvalues of the transfer matrix of the hard square model
  with cylindrical identification using a
  related, but more easily analysed model.
\end{abstract}

\section{Introduction}
\begin{figure}[htbp]
  \begin{center}
    \input{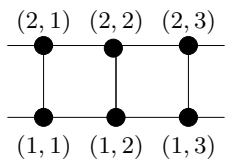}
  \end{center}
  \caption{An undirected graph $G$.}
  \label{fig:intro}
\end{figure}

Let $G = (V,E)$ be an undirected graph.
An \emph{independent set} $I$ in $G$ is a subset of $V(G)$ such that
$u, v \in V(G)$ implies $(u,v) \not\in E(G)$.
Now, take $G$ to be the graph in Figure~\ref{fig:intro}, where
the half-edges are used
to indicate that $(1,1)$ and $(1,3)$ are adjacent and that
$(2,1)$ and $(2,3)$ are adjacent.

We are interested in the independent sets of $G$ and start by listing them.
We begin with the smallest set, the empty set $\emptyset$, which of course
is an independent set in any graph.
Then, the singleton sets $\{ (1,1) \}$, $\{ (1,2) \}$, \ldots, $\{ (2,3) \}$
are all independent.
Finally, for each of the vertices $(1,1), (1,2)$ and $(1,3)$ we can
choose from two vertices of $(2,1), (2,2)$ and $(2,3)$ to form
a total of $6$ independent sets of cardinality 2.
We note that the alternating sum over all independent sets $\sigma$,
\begin{equation} \label{eq:first}
\sum_{ \sigma } (-1)^{|\sigma|} = 1.
\end{equation}
It turns out that this happens a lot for graphs similar to $G$.
Let the vertex set of $C_{m,n}$ be $(i,j)$ for $1 \leq i \leq m$ and
$1 \leq j \leq n$ and let there be an edge between 
$(i_1,j_1)$ and $(i_2,j_2)$ if $i_1 = i_2$ and $|j_1-j_2| = 1$ or
if $|i_1-i_2| = 1$ and $j_1 = j_2$.
Finally, let there be edges between each pair of vertices $(i,1)$ 
and $(i,n)$.
With this construction, we have $G = C_{2,3}$.
It is conjectured \cite{Jonsson06} that when $m$ is odd
and 3 is not a divisor of $\gcd(m-1,n)$, then the sum
$(\ref{eq:first})$ over all independent sets $\sigma$ of $C_{m,n}$ equals 1
and that if $m$ is odd and $3$ does divide $\gcd(m-1,n)$, then
this sum equals -2.
\smallskip

A \emph{hard particle model} in statistical mechanics is 
defined on an underlying graph.
A configuration of the model is a set of particles placed on
the vertices of this graph such that no two particles are
adjacent.
One imagines that the particles have a certain size, or shape, that
extends to their neighbouring vertices in the graph.
Furthermore, the particles may not overlap, as they are considered
\emph{hard}. 
On a square grid, one talks about the hard square model,
because of the shape of the particles.
On a hexagonal grid, the particles are triangles as each vertex
have three neighbours.
The correspondence between configurations of hard particle models
and independent sets on the underlying graph is clear.
\smallskip

The remarkable properties of the hard particle model
on square grid graphs with cylindrical identifications were
first observed by Fendley, Schoutens and van Eerten in \cite{FSvE05}.
They also studied the hard square model with
doubly periodic boundary.
Such identifications have the advantage that they make the
graph vertex-transitive, i.e,
the graph \emph{looks} the same from every vertex.
In this way the influence from the boundary is removed and these
graphs are in that sense closer to the infinite square grid than
the graphs obtained by simply taking an $m$ by $n$ subgraph of
this grid.
Fendley et al. conjectured that for the toroidal identifications,
the roots of the characteristic polynomial of the associated
transfer matrix were all roots of unity.
Furthermore, they conjectured that when $m$ and $n$ were taken
to be relatively prime, the alternating sum $(\ref{eq:first})$ vanished.
These conjectures were proven in 2006 by Jonsson \cite{Jonsson06} 
in a purely combinatorial way.
\smallskip

To determine the partition function of these models
at activity -1 has a direct correspondence in determining
the Euler characteristic of the independence complex of
the underlying graph.
The matchings used in \cite{Jonsson06} were however not suitable for
further topological analysis of the independence complex.
Bousquet-Mélou, Linusson and Nevo \cite{BMLN07} studied a different
region of the square grid which turned out to be more easily
handled using Morse matchings.
They also determined the homotopy type of the independence complex
on a parallelogram in the square grid and related these findings
to certain eigenvalues of the transfer matrix for the
hard square model on $C_{m,n}$.

\smallskip
In this paper, we investigate the independence complexes for
the square and hexagonal grids with cylindrical identifications.
We obtain some results for small sizes by explicitly calculating
the homology groups using a computer software package called
Polymake \cite{Polymake}.
We also construct Morse matchings for some classes of complexes
with small circumference.
The results obtained are enough to make certain conjectures on
the homology groups.
In particular, for odd $n$, there seems to be a symmetry
in the square grid case and we conjecture that for
$i \geq 0$ and $j,k \geq 1$, we have
\[
H_i(I(C_{j,2k+1})) \cong H_i(I(C_{k,2j+1})).
\]
However, it is also clear that the homology does not behave
nicely in many cases, even for odd $n$, which clearly
reduces the applicability of this approach.

The paper is organised as follows.
In Section~\ref{sec:prel} we present some basic notions and results
from algebraic topology, discrete Morse theory and hard particle models.
We use a particular type of Morse matching to determine the homology of
some independence complexes related to 
the hard square and hard triangle models in Section~\ref{sec:explicit}.
In Section~\ref{sec:matching}, we introduce a different matching on the
face poset of these independence complexes, which, while not being
a Morse matching, still turns out to give some information on
the alternating sum of the independence complexes of $C_{m,n}$.
Finally, in Section~\ref{sec:transfermatrices} we review an idea
from \cite{BMLN07} which exploits the knowledge of the partition
function of a hard square model on a different region in order
to try and obtain parts of the spectrum of the transfer matrix
for $C_{m,n}$.

\section{Preliminaries} \label{sec:prel}

In this section we define simplicial complexes, and in particular
independence complexes on grid graphs.
We also review the basic results from discrete Morse theory and
give a framework for constructing Morse matchings on independence
complexes.
Finally, we talk about the connection between hard particle
models from statistical mechanics and independence
complexes on grid graphs.
For a detailed introduction to algebraic topology,
simplicial complexes and homology,
see for example \cite{Munkres} by Munkres.
A classic reference for models in statistical mechanics
is \cite{Baxter82}.
\smallskip

Throughout the paper we will use the following notation for
addition and deletion of a single element $x$ to or from a set $S$.
\[
S + x := S \cup \{x\}, \qquad S - x := S \setminus \{x\}.
\]

To be precise about our base cases, we define the
\emph{Fibonacci numbers} $F_i$ in the following way.
\[
F_0 = 0, F_1 = 1, \text{ and } F_{i+2} = F_{i+1} + F_i \text{ for } i \geq 0.
\]

\subsection{Simplicial Complexes}

Let $X$ and $Y$ be two topological spaces.
If there is a continuous map $f : X \rightarrow Y$,
with a continuous inverse,
then $X$ and $Y$ are said to be \emph{homeomorphic}, denoted by
$X \cong Y$.

Two continuous maps $f : X \rightarrow Y$ and $g : X \rightarrow Y$
are called \emph{homotopic} if there is a continuous map
$F : X \times [0,1] \rightarrow Y$ such that
$F(x,0) = f(x)$ and $F(x,1) = g(x)$.
The spaces $X$ and $Y$ are said to be \emph{homotopy equivalent},
denoted by $X \simeq Y$, if there
are maps $f : X \rightarrow Y$ and $h : Y \rightarrow X$ such
that $f \circ h$ is homotopic to the identity map on $Y$
and $h \circ f$ is homotopic to the identity map on $X$.
A space $X$ that is homotopy equivalent to a point is called \emph{contractible}. This will be denoted by $X \simeq \bullet$.
\smallskip

An {\em (abstract) simplicial complex} $\Delta$, is a collection of subsets over a ground set $V$ such that if $\sigma \in \Delta$ and $\tau \subseteq \sigma$, then $\tau \in \Delta$.
The elements of $\Delta$ are called {\em faces}, and the maximal faces are called {\em facets}.
The {\em dimension} of a face $\sigma$ is defined to be $\dim (\sigma) = |\sigma|-1$.

The {\em reduced Euler characteristic} ${\tilde \chi}$ of a simplicial complex $\Delta$ is the following sum.
\begin{equation} \label{eq:euler}
{\tilde \chi}(\Delta) := \sum_{\sigma \in \Delta} (-1)^{\dim (\sigma)}
\end{equation}

\medskip
Given a simplicial complex $\Delta$, the {\em homology group} $H_i(\Delta)$
of $\Delta$ in dimension $i$, over a ring $R$, is defined by 
$H_i(\Delta) = C_i(\Delta)/B_{i}(\Delta)$.
Here, the groups $C_i$ and $B_{i}$ are the groups of 
$i$-chains and $i$-boundaries, respectively.
The rank of $H_i(\Delta)$ is called the $i$th \emph{Betti number}
and is denoted by $\beta_i$.
The following relation to the Euler characteristic, which can
alternatively be taken as its definition, 
is called the \emph{Euler-Poincaré formula}.
\begin{equation} \label{eq:eulerpoincare}
{\tilde \chi}(\Delta) = \sum_{i \geq -1} (-1)^i \beta_i
\end{equation}
Here, $\beta_{-1}$ is defined to be $1$ if the complex
$\Delta = \{ \emptyset \}$ and $0$ otherwise.
Homotopy equivalent spaces have isomorphic homology groups
so ${\tilde \chi}$ is a topological invariant.

For a topological space $X$, we let $\partial X$ denote the
\emph{boundary} of $X$ and $\text{int} \ X = X \setminus \partial X$
the \emph{interior} of $X$.
The $n$-dimensional simplex $\Delta_n$ is the full simplicial complex
on a $n+1$-point ground set, i.e, $\Delta_n = 2^V$, with $|V| = n+1$.
The $n$-dimensional \emph{ball} $B^n$ is the subspace of the
$n$-dimensional Euclidean space ${\mathbb R}^n$ given by
$B^n := \{ \vv \in \mathbb{R}^n \;|\; \|\vv\| \leq 1 \}$.
The $n$-dimensional \emph{sphere} $S^n$ is similarly given by
$S^n := \{ \vv \in \mathbb{R}^{n+1} \;|\; \|\vv\| = 1 \}$.
That is, $S^n = \partial B^{n+1}$.
Apparently, we have $B^n \cong \Delta_n$ and 
$S^n \cong \partial \Delta_{n+1}$,
where the simplicial complex $\partial \Delta_{n+1}$ 
can be described by $2^V \setminus V$, for some $|V| = n+1$.

\medskip
Let $\Delta$ be a simplicial complex.
The {\em cone} over $\Delta$, denoted by $\mathrm{cone}(\Delta)$, 
is the simplicial complex over the ground set of $\Delta$ together with a new vertex $c$, the cone point.
\[
\mathrm{cone}(\Delta) = \{ \sigma, \sigma + c \;|\; \sigma \in \Delta\}
\]
The cone over any simplicial complex is contractible.
\medskip

The {\em suspension} over $\Delta$, denoted by $\mathrm{susp}(\Delta)$,
is the following simplicial complex over the ground set of $\Delta$ together with two new vertices $s_0$ and $s_1$.
\[
\mathrm{susp}(\Delta) = \{ \sigma, \sigma + s_0, \sigma + s_1  \;|\; \sigma \in \Delta\}
\]
For any simplicial complex $\Delta$ and $i \geq 0$,
\begin{equation} \label{eq:susp}
  H_i(\Delta) \cong H_{i+1}\left(\mathrm{susp}(\Delta)\right).
\end{equation}

\subsection{Independence Complexes of Regular Grids} \label{subsec:reggrid}

\begin{figure}[htbp]
  \begin{center}
    \input{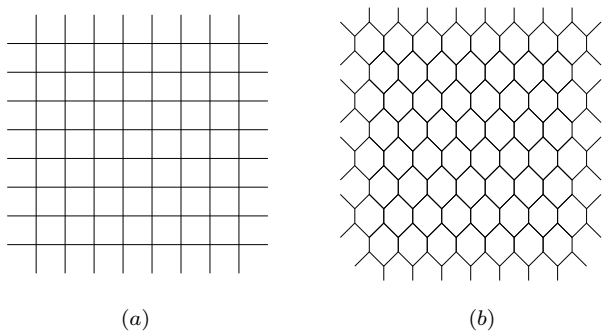}
  \end{center}
  \caption{(a) The square grid. (b) The hexagonal grid.}
  \label{figure:grids}
\end{figure}

Regular grids are tilings of the 
2-dimensional plane by regular convex polygons.
There exist three types of regular tilings, triangular, square and
hexagonal.
In this paper we have considered the square and the hexagonal grid.
These are depicted in Figure~\ref{figure:grids}.

Let $S$ be the graph induced by the boundary of the square grid.
Then, $S$ is isomorphic to an infinite graph on the vertex set 
${\mathbb Z} \times {\mathbb Z}$,
with edges between $(u_1,u_2)$ and $(v_1,v_2)$ if $u_1 = v_1$ and $|u_2-v_2| = 1$ or
if $|u_1-v_1| = 1$ and $u_2 = v_2$.
From now on, we will assume that $S$ is defined in this way.
Let $S_{m,n}$ be the graph induced from $S$ by the vertex set $[m] \times [n]$,
where $[m]$ denotes the set $\{1,\ldots,m\}$.
By identifying sides of $S_{m,n}$, we can obtain graphs embeddable on different
2-dimensional topological surfaces.
Specifically, let $C_{m,n}$ be $S_{m,{n+1}}$, with vertices $(i,n+1)$ and $(i,1)$ identified, for $i \in [m]$.
Furthermore, let $T_{m,n}$ as the graph $S_{m+1,n+1}$, with the vertices
$(i,n+1)$ and $(i,1)$ as well as the vertices $(m+1,j)$ and $(1,j)$ identified,
for $i \in [m]$ and $j \in [n]$, respectively.
The graph $C_{m,n}$ is obviously embeddable on a cylinder
(but not on the plane), while $T_{m,n}$
is embeddable on a torus (but not on a cylinder).

\begin{figure}[htbp]
  \begin{center}
    \input{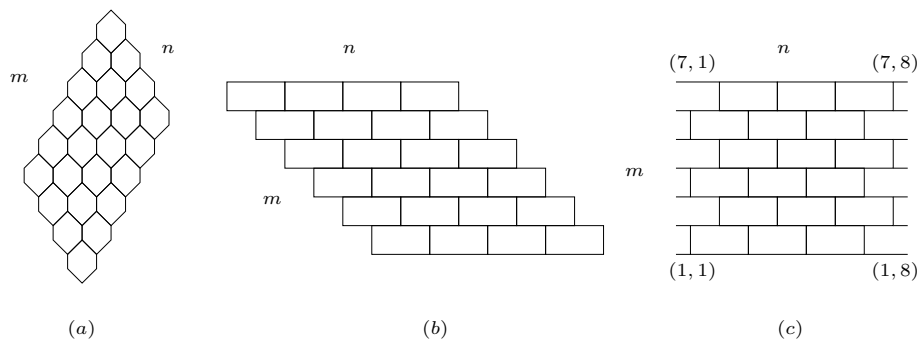}
    \caption{$(a)$ $H_{6,4}$ as a subgraph of $H$. $(b)$ $H_{6,4}$ as a brick wall. $(c)$ The cylinder $C^H_{6,4}$.}
    \label{figure:hexgrid}
  \end{center}
\end{figure}

\medskip
Let $H$ be the graph induced by the boundary of the hexagonal grid.
Similarly to the case of $S$, we can find an appropriate representation of
the graph $H$ with vertex set ${\mathbb Z} \times {\mathbb Z}$.
This graph has edges between $(u_1,u_2)$ and $(v_1,v_2)$ if
$u_1 = v_1$ and $|u_2-v_2| = 1$ or if $u_1-v_1 = 1$, $u_2 = v_2$,
and $u_1 + u_2 \equiv_2 0$.
For positive integers $m$ and $n$ let $H_{m,n}$ be the subgraph of $H$
shown in Figure~\ref{figure:hexgrid}$(a)$.
We can straighten this graph and draw it on the square grid, as
is shown in $(b)$.
We will use this \emph{brick wall}-like appearance in this paper.
If we identify the left and right boundaries of the graph in
Figure~\ref{figure:hexgrid}$(b)$, we get a graph $C^H_{m,n}$
with the topology of a cylinder.
Figure~\ref{figure:hexgrid}$(c)$ shows an alternative way of
drawing this cylinder, using half edges on the left and
right boundaries to indicate the identifications.
We will use the vertex set 
$[m] \times [n] \subseteq {\mathbb Z} \times {\mathbb Z}$
to represent these graphs, with the vertices numbered
as in Figure~\ref{figure:hexgrid}$(c)$.
We will also briefly look at the hexagonal grid with toroidal
identifications.
In this case, we will use the representation in Figure~\ref{figure:hexgrid}$(b)$, identifying the left and right boundaries as well as the top and bottom ones.
We denote the graph obtained by these identifications by $C^T_{m,n}$.
Note that this graph is different from the one we would obtain by doing
the top to bottom identifications in Figure~\ref{figure:hexgrid}$(c)$.

\smallskip
We will be studying simplicial complexes related to these regular grids.
Let $G = (V, E)$ be a finite, simple and undirected graph.
An \emph{independent set} $I$ in $G$ is a subset of the vertex set
$V(G)$ such that if $u, v \in I$ then $(u,v) \not\in E(G)$.
That is, no pair of neighbours may be in $I$ at the same time.

Define the \emph{independence complex} $I(G)$ as
the simplicial complex on the ground set
$V(G)$ with $\sigma \in I(G)$ if and only if
$\sigma$ is an independent set in $G$.
Similarly,
let $X(G)$ be the \emph{flag complex} of $G$.
This is the simplicial complex on the ground set
$V(G)$ with $\sigma \in X(G)$ if and only if
$\sigma$ induces a complete subgraph in $G$.
It is obvious from the definitions that $I(G) = X({\bar G})$,
where ${\bar G} = (V, {\bar E})$ denotes the complement
graph of $G$ with respect to the complete graph on $V$, i.e, ${\bar E} = \binom{V}{2} \setminus E$.

In general, computing $I(G)$ is believed to be hard, since
the facets of $I(G)$ are the maximal independent sets of $G$.
Any such facet of maximal dimension therefore gives a solution
to the {\em maximum independent set problem}, which is a well known
{\sc NP}-complete problem \cite{Karp72}.
Dual to this is the {\sc NP}-complete {\em maximum clique problem} of finding
a complete subgraph (or \emph{clique}) in $G$ of maximal size.
These cliques are given by the highest dimensional facets of $X(G)$.


\subsection{Discrete Morse Theory}

\emph{Discrete Morse theory} was introduced in \cite{Forman95} by Robin Forman.
It relates the homotopy type of a simplicial complex to that of a
\emph{CW-complex} with (ideally) simpler structure.
See \cite{Forman02} for a hands-on introduction to this subject.

Let $\Delta$ be a simplicial complex.
We define the \emph{face poset} $P(\Delta)$ on $\Delta$ as
the set of faces in $\Delta$ ordered by inclusion.
The \emph{underlying graph} of the Hasse diagram of $P(\Delta)$,
or simply the underlying graph of $P(\Delta)$,
is the directed graph with  the faces of $\Delta$ as vertices and
the cover relations of $P(\Delta)$ as edges, directed from the
smaller faces to the larger faces.
Let $M$ be a matching on the underlying graph of $P(\Delta)$.
We will also say that $M$ is a matching on the complex $\Delta$.
The \emph{graph induced by $M$} is the graph obtained 
from the underlying graph of the Hasse diagram by reversing
the direction of the edges which are in $M$.

\begin{definition}
  A \emph{Morse matching} $M$ on a simplicial complex $\Delta$ is a (partial)
  matching of the Hasse diagram of $P(\Delta)$ such that the graph
  induced by $M$ is acyclic.
\end{definition}

The main theorem of Morse theory relates the homotopy type of a
simplicial complex to that of a certain \emph{CW-complex} related
to the unmatched faces of a Morse matching.
These unmatched faces are sometimes referred to as the \emph{critical cells}
of the matching.
For our purposes, the following definition of a finite CW-complex will
suffice.
\begin{definition}
  A finite CW-complex (or cell complex) $C$ is a collection of cells
  such that for each cell $e_\alpha$
  \begin{itemize}
    \item[(a)] There is a homeomorphism from $B^k$ to $e_\alpha$ for
      some $k$.
    \item[(b)] $e_\alpha \setminus \text{int } e_\alpha$ is the union of cells
      $e_\beta$ in $C$.
  \end{itemize}
\end{definition}

The main theorem of discrete Morse theory can be stated as follows,
using the notion of a Morse matching.
\begin{theorem} \label{thm:morse}
  Let $\Delta$ be a simplicial complex with a Morse matching $M$.
  Assume that for each $p \geq 0$, there are $u_p$ unmatched p-simplices.
  Then $\Delta$ is homotopy equivalent to a CW-complex with exactly
  $u_p$ cells of each dimension $p$, for each $p > 0$, and
  $u_0+1$ cells of dimension $0$.
\end{theorem}

\begin{cor} \label{cor:wedgeofspheres}
  Let $\Delta$ be a simplicial complex with a Morse matching $M$
  such that $u_p = 0$ for all but one $p$.
  Then, for this particular $p$,
  \[
  \Delta \simeq \bigvee_{u_p} S^p.
  \]
  That is, $\Delta$ is homotopy equivalent to a wedge of $u_p$
  $p$-dimensional spheres.
\end{cor}

\begin{cor} \label{cor:perfectmorse}
  Let $\Delta$ be a simplicial complex with a perfect Morse matching $M$.
  Then, $\Delta$ is contractible.
\end{cor}

In this paper, we will apply discrete Morse theory to
independence complexes.
The acyclic matchings that we will use are all instances of the
following general construction.
By describing them in this framework, we can easily determine
their acyclicity.
The construction is similar to the \emph{matching trees}
of \cite{BMLN07}.
Throughout this paper, we will however use the term
\emph{matching tree} to refer to our construction.

\begin{example} \label{ex:mtree}
  The following example shows the basic idea behind the matching trees.
  In Figure~\ref{fig:mtree1}, we construct a matching on the complex
  $I(S_{3,2})$.
  First we partition the set of faces in to two disjoint sets, based on
  whether the vertex $(1,1)$ is in the face or not.
  We indicate this by writing $(1,1)$ to the left of the root in the tree.
  To the right of the root we place the graph $S_{3,2}$ which indicates
  that all faces of the complex are still unmatched.
  In the left subtree we have forced $(1,1)$ not to appear in any face.
  We illustrate this by placing an empty particle on $(1,1)$ in the
  graph to the right of the root of this subtree.
  In the right subtree we have instead forced $(1,1)$ to be present
  in every face.
  This is illustrated by a filled particle on $(1,1)$ and empty
  particles on its neighbours $(1,2)$ and $(2,1)$.
  In both subtrees, this effectively reduces the problem to the
  construction of a matching on two smaller graphs,
  induced from $S_{3,2}$ by the complement of the fixed vertices.
  In the next step, the vertex $(1,2)$ in the left graph, and $(2,2)$
  in the right graph both have degree 1.
  We use these vertices to create the matchings $M(1,2)$
  and $M(2,2)$ in the respective subtree.
  Let $\Sigma_l = I(S_{3,2} - (1,1))$ and $\Sigma_r = I(S_{3,2} \setminus \{ (1,1), (1,2), (2,1) \})$, and define the matchings as follows.
  \[
  M(1,2) := \{(\sigma, \sigma - (1,2)) \;|\; (1,2) \in \sigma \in \Sigma_l \},
  \]
  \[
  M(2,2) := \{(\sigma, \sigma - (2,2))  \;|\; (2,2) \in \sigma \in \Sigma_r \}.
  \]
  The unmatched faces $\sigma$ remaining after removing the faces matched by
  $M(1,2)$ must have $(1,2), (2,1), (3,2) \not\in \sigma$ 
  and $(2,2) \in \sigma$.
  This further reduces the graph.
  In the right subtree, it completely fixes one single graph.
  In the left subtree, there is still a vertex $(3,1)$ on which we can match.
  Matching $\{(2,2)\}$ with $\{(2,2),(3,1)\}$ 
  removes all the faces from this subtree, which we indicate by writing
  $\emptyset$ to the right of the last node.
  \begin{figure}[htbp]
    \begin{center}
      \input{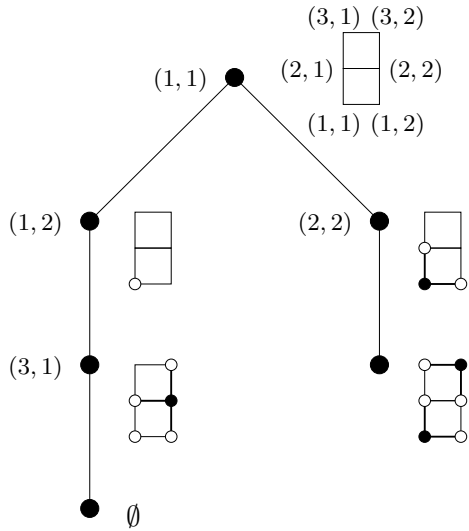}
      \caption{A matching tree for the simplicial complex $I(S_{3,2})$.}
      \label{fig:mtree1}
    \end{center}
  \end{figure}

  In Proposition~\ref{prop:tree}, we will see that the union of the
  matchings made in this type of construction forms an acyclic matching
  on the original complex.
  We can therefore apply Theorem~\ref{thm:morse}, and due to the
  remaining 1-dimensional face in the right subtree, which is
  the only critical cell, we conclude
  that $I(S_{3,2})$ is homotopy equivalent to a 1-dimensional
  sphere.
\end{example}

\begin{definition}[Matching tree]
Let $T$ be a finite, rooted tree in which each internal node
has either one or two children.
To each internal node (which includes the root)
$t \in T$, we associate a vertex $p(t) \in G$.

We say that $T$ is a \emph{matching tree} if
the set of vertices $p(t)$ fulfil the following two conditions.
For any node $t_r \in T$, let $(t_1, \ldots, t_r)$ be
the path from the root $t_1$ to $t_r$.
Then,
\begin{itemize}
\item all $p(t_1), \ldots, p(t_r)$ are distinct, and
\item $\{ p(t) \;|\; t \in \{t_1, \ldots, t_r\} \cap V_m\}$ forms an independent set in $G$,
\end{itemize}
where $V_m$ denotes the set of nodes with one child.
\end{definition}

The vertices $p(t)$ are called \emph{the pivots} in \cite{BMLN07}.
They are the vertices that we wrote to the left of the nodes in
Figure~\ref{fig:mtree1}.
Given the complex $I(G)$ and the tree $T$, we determine a matching
on $P(I(G))$ as follows.
To each node $t \in T$, we associate a \emph{content} $C(t)$, which is
a subset of the complex $I(G)$.
This is the set that was illustrated by the graphs to the right
of the nodes in Figure~\ref{fig:mtree1}.
For the root $r \in T$, we let $C(r) = I(G)$, and for the other nodes
the content will be determined recursively.

\smallskip
If $t$ has two children, then we call $t$ a \emph{split node}, and we let
the content of the root of the left subtree of $t$ be the set of faces
$\{ \sigma \in C(t) \;|\; p(t) \not\in \sigma\}$ and the content of the
root of the right subtree of $t$ be the faces
$\{ \sigma \in C(t) \;|\; p(t) \in \sigma\}$.

\smallskip
If $t$ has only one child, then we call $t$ a \emph{match node}.
We define the partial matching $M(t)$ associated to the match node
and the set of matched faces $F(t)$ as follows.
\[
M(t) := \{ (\sigma, \sigma - p(t)) \;|\; p(t) \in \sigma \in C(t) \},
\qquad
F(t) := \bigcup_{(\sigma,\tau) \in M(t)} \{\sigma, \tau\}
\]
We will see below that the set $M(t)$ 
forms a partial matching on $C(t)$.
For match nodes,
we let the content of the child of $t$ be $C(t) \setminus F(t)$.

\smallskip
\begin{definition}
  Let $G$ be a graph and $\sigma \in I(G)$.
  We say that a position $x \in V(G)$ is \emph{free in $\sigma$}
  if $x \not\in \sigma$ and $\sigma + x$ is an independent set in $G$.
\end{definition}
The following lemma provides all the structure we need on the
sets $C(t)$ to show acyclicity for our matchings.

\begin{lemma2} \label{lemma:hake}
  Let $T$ be a matching tree and let $t \in V_m$.
  If $\sigma \in C(t)$ and $p(t)$ is free in $\sigma$, 
  then $\sigma + p(t) \in C(t)$
  Conversely, if $\sigma \in C(t)$ and $p(t) \in \sigma$, then
  $\sigma - p(t) \in C(t)$.
\end{lemma2}

\begin{proof}
  Both parts of the lemma obviously hold true for the root node.
  We will show that if either part fails for some $t$, 
  where $t$ is taken
  to be the first node on the path from the root to $t$ for which this
  is the case, then
  the other part must fail at an earlier node $t'$, closer to the root.
  If we apply this argument twice, we reach a node
  $t''$ which is strictly between the root and $t$ and
  for which the same part of the lemma fails.
  This contradicts our choice of $t$, and the lemma follows.

  \smallskip
  Assume that $\sigma \in C(t)$ and that $p(t)$ is free in $\sigma$.
  The reason for $\sigma + p(t)$ not to be in $C(t)$ could be
  that there was a split node which had $\sigma$ in one subtree
  and $\sigma + p(t)$ in the other.
  This can not happen, however, since the two faces only differ
  in $p(t)$ and along the path from the root, all pivots must
  be distinct.

  The only remaining possibility is then that $\sigma + p(t)$ was
  matched in an earlier node $t'$.
  Now, if $p(t')$ is free in $\sigma + p(t)$, then it is
  also free in $\sigma$, and
  by our assumption, $\sigma$ would have been matched at $t'$ as well.
  So, we can assume that $p(t') \in \sigma + p(t)$, and then
  also $p(t') \in \sigma$.
  The second part of the lemma therefore fails at $t'$.

  \smallskip
  For the second part, assume that $\sigma \in C(t)$ and $p(t) \in \sigma$,
  but that $\sigma - p(t) \not\in C(t)$.
  Again, $\sigma - p(t)$ could not have been removed by an earlier
  split node since it only differs in $p(t)$ from $\sigma$.

  So assume that $\sigma - p(t)$ was matched in an earlier node $t'$.
  If $p(t') \in \sigma - p(t)$, then $p(t') \in \sigma$ and by
  induction, $\sigma$ would have been matched with $\sigma - p(t')$
  at $t'$.
  Therefore, $p(t')$ is free in $\sigma - p(t)$.
  By the restriction that the pivots of the match nodes form an independent
  set in $G$, we have that $p(t')$ is free in $\sigma$ as well.
  The first part of the lemma therefore fails at $t'$.
  \qed
\end{proof}

Using Lemma~\ref{lemma:hake}, it is easy to see that in fact, 
for $t \in V_m$, we have
\[
F(t) \subseteq C(t).
\]
That is, for each match node $t \in V_m$, $M(t)$ is a partial matching
on the subgraph of the underlying graph of $P(I(G))$
induced by the set $C(t)$.
It is acyclic, as the only way to ``move down'' in the 
directed graph induced by $M(t)$ is by removing $p(t)$.

It therefore makes sense to define the matching $M$ determined by the
tree $T$ as
\begin{equation} \label{eq:M}
M := \bigcup_{t \in V_m} M(t).
\end{equation}
We will now show that $M$ is acyclic, so that
Theorem~\ref{thm:morse} is applicable.
The critical cells will be those in the set
$\bigcup_{t \in V_l} C(t)$,
where $V_l$ denotes the set of leaves in $T$.
Our proof will use the following lemma, 
which appears in \cite{Jonsson:thesis}.
\begin{lemma2}[Cluster Lemma] \label{lemma:cluster}
  Let $\Delta$ be a simplicial complex and $f : P(\Delta) \rightarrow Q$ be
  a poset map, where $Q$ is an arbitrary poset.
  For $q \in Q$, let $M_q$ be an acyclic matching on $f^{-1}(q)$.
  Let
  \[
  M = \bigcup_{q\in Q} M_q.
  \]
  Then, $M$ is an acyclic matching on $\Delta$.
\end{lemma2}



The following proposition justifies our definition of a matching tree,
showing that any tree gives an acyclic matching on the corresponding
independence complex.
The proof is entirely guided by the Cluster lemma (Lemma~\ref{lemma:cluster}).
\begin{proposition2} \label{prop:tree}
  Let $G$ be a graph and $T$ a matching tree defined on
  the independence complex $I(G)$.
  Then, the matching $M$ given by $(\ref{eq:M})$
  on $I(G)$ is acyclic.
\end{proposition2}

\begin{proof}
  We start by defining a total order on the nodes of the tree $T$.
  First, if $t'$ is an ancestor of $t$, then let $t \leq_T t'$.
  Then, if $T_1$ and $T_2$ are two subtrees of a node $t$,
  with $T_1$ being the subtree where $p(t) \not\in \sigma$ for
  any $\sigma \in C_(t_1)$,
  let $t_1 \leq_T t_2$ for all $t_1 \in T_1$ and $t_2 \in T_2$.
  
  Let $Q = (V(T), \leq_T)$ and 
  \[
  f^{-1}(t) = \begin{cases}
    F(t) & \text{if $t \in V_m$,} \\
    C(t) & \text{if $t \in V_l$,} \\
    \emptyset & \text{otherwise.} \\
  \end{cases}
  \]
  It is easy to check that the $F(t)$ are pairwise disjoint and that
  $F(t_1) \cap C(t_2) = \emptyset$ for all $t_1 \in V_m$ and $t_2 \in V_l$.
  Furthermore, each face of $I(G)$ is either in $F(t_1)$ for some
  $t_1 \in V_m$ or in $C(t_2)$ for some $t_2 \in V_l$.
  Thus, $f : I(G) \rightarrow Q$ is a uniquely defined map.
  It remains to verify that it is indeed a poset map.
  
  Let $\sigma, \tau \in I(G)$, with $\sigma \subseteq \tau$.
  If $f(\sigma)$ and $f(\tau)$ are in different subtrees of some node $t$,
  then $\sigma$ must be in the left, as we know that $p(t)$ is in all faces
  of the right subtree, but in none of the left.
  Thus, $f(\sigma) \leq_T f(\tau)$.

  If instead $f(\sigma)$ and $f(\tau)$ are on the same path from the
  root to some leaf, then we will assume that $f(\tau)$ is closer to the
  root than $f(\sigma)$.
  Let $v = p(f(\tau))$ be the pivot at $f(\tau)$.
  If $v \in \sigma$, then Lemma~\ref{lemma:hake} says that
  $\sigma$ should have been matched in the node $f(\tau)$,
  so we have a contradiction.
  Therefore, $v \not\in \sigma$.
  But then $\sigma \subseteq \tau-v$, and $\tau-v$ was matched with
  $\tau$.
  This implies that $v$ is free in $\sigma$, and Lemma~\ref{lemma:hake}
  again shows that $\sigma$ should have been matched in $f(\tau)$.
  We have a contradiction. Therefore $f(\tau)$ must be
  equal to or lie further from the root than $f(\sigma)$, 
  which means that $f(\sigma) \leq_T f(\tau)$.

  Finally, for each $t \in V_m$, we have that 
  $M(t)$ is a perfect, acyclic matching of
  $f^{-1}(t) = F(t)$.
  We can therefore use Lemma~\ref{lemma:cluster} to
  conclude that $M$ is acyclic on the entire $I(G)$.
  \qed
\end{proof}  

We finally note that our construction disallows some strategies
for choosing split and match nodes which would have resulted in acyclic
matchings.
This is due to the constraint that the match nodes in each
path from the root must form an independent set.
This condition is not strictly necessary, but
for the purposes of this article, this construction will suffice.

\subsection{Hard particle models} \label{subsec:hardp}

A particular type of model in two-dimensional statistical mechanics
is the {\em hard particle model}.
In such a model, given a regular grid like the ones described in 
Section~\ref{subsec:reggrid}, a configuration is a set of particles
placed on the sites of the grid so that no two particles are adjacent.
The name is explained by replacing each particle with a polygon,
using the neighbours of the original particle as the vertices of the
polygon.
This polygon is assumed to be ``hard'', meaning that no two
polygons may overlap in a configuration.
When the grid is the square grid, the polygons are squares, and the
model referred to as the {\em hard square model}.
For the hexagonal grid, the polygons are triangles, and one has the
{\em hard triangle model}.
For a thorough introduction to models in statistical mechanics,
see \cite{Baxter82}, in which can also be found the solution to the
hard particle model on a triangular grid.
This solution was first presented by Baxter in \cite{Baxter80} in 1980.

To a model, one can associate a {\em partition function} $Z$.
This is the generating function of the configurations with 
respect to the number of particles of a configuration.
\begin{equation} \label{eq:pf}
Z = \sum_{\sigma} z^{|\sigma|},
\end{equation}
where the sum is taken over all allowed configurations of the model.
The variable $z$ of $Z$ is called the {\em activity} and controls the
behaviour of the model.


\smallskip

A commonly used framework for deriving the partition function is the
\emph{transfer matrix method}.
See \cite[Chapter 4.7]{Stanley:EnCo1} for an introduction to this method.
Here, we illustrate the idea by giving a transfer matrix
for the rectangle $S_{m,n}$.
More examples are given in Section~\ref{sec:transfermatrices}.

Let $L_n$ denote the set of independent sets on the line $S_{1,n}$.
Since $|L_1| = 2$, $|L_2| = 3$, and the recursive relation
$|L_{n+2}| = |L_{n+1}| + |L_{n}|$ holds, 
we have $|L_n| = F_{n+2}$, the Fibonacci numbers.
Let $T(z) := T(z)_{\sigma \tau}$ be the $F_{n+2} \times F_{n+2}$ matrix
where the individual matrix elements are given by
\[
T(z)_{\sigma \tau} = \begin{cases} z^{|\tau|} & 
  \text{if $\sigma \cap \tau = \emptyset$} \\
  0 & \text{otherwise.} \end{cases}
\]
The matrix $T(z)$ is called a transfer matrix of the model.
It describes the effect of adding an additional row.
We will only be studying models at activity -1
so we leave out the variable $z$ and write $T = T(-1)$.
The generating function of $Z(S_{m,n})$ with respect to $n$ can
now be written as
\[
\sum_{n \geq 0} Z(S_{m,n}) t^n = (1-tT)^{-1}.
\]
In addition, by simply taking the trace of the $n$th power of the matrix $T$,
we can also express the partition function $Z(C_{m,n})$ of a cylinder.
\[
Z(C_{m,n}) = \mathrm{tr}(T(z)^n),
\]
or as a generating function
\[
  \sum_{n \geq 0} Z(C_{m,n}) t^n = \mathrm{tr}(1-tT)^{-1}.
\]
The \emph{characteristic polynomial} $P(t)$ of $T$ is given by
\[
P(t) = \det(tI-T).
\]
The roots of the characteristic polynomial can also be found by taking the poles of the generating
function $\mathrm{tr}(1-tT)^{-1}$.

It was observed by Fendley, Schoutens and van Eerten \cite{FSvE05}
that the hard square model (on a torus) exhibit some interesting properties 
when the activity is chosen to be $-1$.
Among other things, they conjectured that all the roots of the
characteristic polynomial are roots of unity.
These conjectures were rigorously proven by Jonsson \cite{Jonsson06}
by means of applying an ingeniously designed series of matchings on the set
of configurations.
Doing this, Jonsson discovered an underlying structure of a particular type
of rhombus tilings of the plane which he used to settle the conjectures.
An attempt to mimic the idea of Jonsson is carried out in Section~\ref{sec:matching},
albeit with limited results.

\smallskip
The following conjecture for the hard particle model on a cylinder
in the square grid at activity -1 can be found in \cite{Jonsson06}.
The conjecture is based on results for $m \leq 11$.
\begin{conjecture} \label{conj:cylinder}
  For odd $n$,
  \[
  Z(C_{m,n}) = \begin{cases}
    -2 & \text{if $3\,| \gcd(m-1,n)$,} \\
    1 & \text{otherwise.} \\
  \end{cases}
  \]
\end{conjecture}

\smallskip
From (\ref{eq:euler}) and (\ref{eq:pf}) we see that, at $z = -1$, 
\[
Z(G) = -{\tilde \chi}(I(G)).
\]
That is, the partition function at activity $-1$ is (up to a negative sign)
equal to the 
reduced Euler characteristic of the
independence complex of the underlying graph.
For this reason, we will sometimes refer to $Z(G)$ as the
\emph{alternating sum} of the corresponding independence complex.
In the next section, we investigate if this connection to topology
can be taken further than just to the reduced Euler characteristic.

\section{Explicit results} \label{sec:explicit}

In this section we present the explicitly calculated homology groups
of some independence complexes on graphs related to the square and the
hexagonal grid.
First, we do some small cases on the square grid by hand, 
for a fixed circumference $n$.
Then, we present some larger cases, where we have used
a computer to calculate the homology groups.
The calculations were made using the software package Polymake \cite{Polymake}
version 2.2 by Michael Joswig and Ewgenij Gawrilow with many contributors.

\subsection{Results for small circumference} \label{subsec:small}

For some cases with small $n$, we can find suitable matching trees
and work out the homotopy type.
We do this here for the square grid when $n \leq 5$ and for
the hexagonal grid when $n = 1$ or $n = 2$.
For the case $n = 4$, we are only able to determine the alternating sum.
This result is given in Proposition~\ref{prop:hexcyl4}.
We note that for the square grid, these results are consistent
with Conjecture~\ref{conj:cylinder}.

\medskip
When $n = 1$ there can be no particles on the grid since every
vertex in the grid is a neighbour to itself.
Therefore $I(C_{m,1}) \homeo \{ \emptyset \}$ and 
$Z(C_{m,1}) = 1$.

\begin{proposition2} \label{prop:sq2}
Let $m \geq 1$. Then, $I(C_{m,2}) \homot S^{\lceil m/2 \rceil-1}$
and $Z(C_{m,2}) = (-1)^{\lceil m/2 \rceil-1}$.
\end{proposition2}

\begin{proof}
Note that $C_{m,2} = S_{m,2}$.
The construction in Example~\ref{ex:mtree} in Section~\ref{sec:prel}
is easily extended to arbitrary $m$.
The only unmatched configuration has
$m/2$ particles when $m$ is even and $(m+1)/2$ particles when $m$ is
odd.
\qed
\end{proof}

\begin{proposition2} \label{prop:sq3}
  Let $m \geq 1$. Then,
  \[
  I(C_{m,3}) \simeq
  \begin{cases}
    S^{2m/3-1} & \text{if $m \equiv_3 0$,} \\
    \bigvee_2 S^{2(m-1)/3} & \text{if $m \equiv_3 1$,} \\
    S^{2(m-2)/3+1} & \text{if $m \equiv_3 2$,} \\
  \end{cases}
  \]
  and
  \[
  Z(C_{m,3}) =
  \begin{cases}
    1 & \text{if $m \equiv_3 0$,} \\
    -2 & \text{if $m \equiv_3 1$,} \\
    1 & \text{if $m \equiv_3 2$.} \\
  \end{cases}
  \]
\end{proposition2}

\begin{figure}[htbp]
  \begin{center}
    \input{mtree5.pstex_t}
    \caption{Matching tree for $C_{m,3} - (1,1)$.}
  \end{center}
  \label{figure:n3}
\end{figure}

\begin{proof}
Let $C'_{m,n} = C_{m,n} - (1,1)$.
We construct an acyclic matching on $C_{m,3}$ as follows.
If we start by splitting on $(1,1)$, we are left with
$C'_{m,3}$ in the left subtree, and after a matching, 
a graph isomorphic to $C'_{m-1,3}$ in the right subtree.
Figure~\ref{figure:n3} shows the matching tree construction for $C'_{m,3}$.
Note that the last remaining graph in the left subtree is
isomorphic to $C'_{m-3,3}$ and that we have added two particles
on the path from the root to this node.

When $m = 1$, we can construct a matching leaving two configurations
with one particle each,
when $m = 2$, we find a perfect matching, and
when $m = 3$, we are left with one configuration, having two particles.
The proposition follows by 
adding the configurations from the two subtrees, and
using Proposition~\ref{prop:tree} and Corollary~\ref{cor:wedgeofspheres}.
\qed
\end{proof}

\begin{proposition2} \label{prop:sq4}
  Let $m \geq 1$. Then,
  \[
  I(C_{m,4}) \simeq
  \begin{cases}
    \bigvee_{m+1} S^{m-1} & \text{if $m$ is even,} \\
    \bigvee_{m} S^{m-1} & \text{if $m$ is odd,} \\
  \end{cases}
  \]
  and
  \[
  Z(C_{m,4}) = 
  \begin{cases}
    m+1 & \text{if $m$ is even,} \\
    -m & \text{if $m$ is odd.} \\
  \end{cases}
  \]
\end{proposition2}

\begin{figure}[htbp]
  \begin{center}
    \input{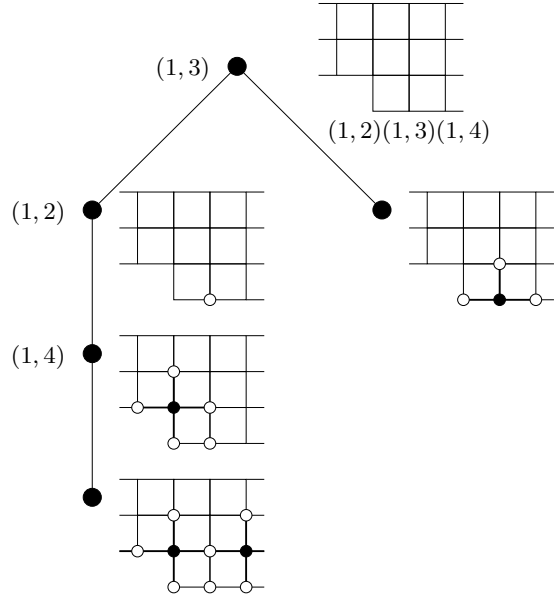}
    \caption{Matching tree for $C_{m,4} - (1,1)$.}
  \end{center}
  \label{figure:n4}
\end{figure}

\begin{proof}
Let $C'_{m,n} = C_{m,n} - (1,1)$.
Figure~\ref{figure:n4} shows how to construct a matching tree for the
graph $C'_{m,4}$.
Note that the left subtree gives a recursive relation of length 2.
When $m$ is even, there will be one configuration with $m$
particles left at the end of this subtree.
When $m$ is odd, we have a perfect matching.
In the right subtree, there remains a graph which is isomorphic to
$C'_{m-1,n}$.
When $m = 1$, this subtree has one configuration with one particle.
Thus, the entire tree leaves 
$\lfloor m/2 \rfloor + 1$ configurations with $m$ particles.

Now, we build the matching tree for $C_{m,4}$ as follows.
We first split on $(1,1)$. The left subtree will be the one
previously constructed for $C'_{m,4}$ with $\lfloor m/2 \rfloor + 1$
remaining configurations.
In the right subtree, we can match on $(3,1)$, leaving
us with a graph isomorphic to $C'_{m-2,4}$, and
$\lfloor (m-2)/2 \rfloor + 1$ remaining configurations.

We then use Proposition~\ref{prop:tree} to show that we have
an acyclic matching with $2 \lfloor m/2 \rfloor + 1$
unmatched configurations, each with $m$ particles.
Corollary~\ref{cor:wedgeofspheres} finally gives us the homotopy type.
\qed
\end{proof}

\begin{proposition2}
  Let $m \geq 1$. Then,
  \[
  I(C_{m,5}) \simeq \begin{cases} 
    S^{m-1} & \text{if $m$ is even,} \\
    S^m & \text{if $m$ is odd,} \\
  \end{cases}
  \]
  and
  \[
  Z(C_{m,5}) = 1.
  \]
\end{proposition2}


\begin{proof}
  We give a sketch of the matching tree for $C_{m,5}$.
  First, we split on both $u = (1,1)$ and $v = (3,1)$.
  For $S \subseteq 2^{\{u,v\}}$, let
  \[
  X_m(S) = \{ \sigma \in I(C_{m,5}) \;|\; \sigma \cap \{u,v\} = S \}.
  \]
  Then, $X_m(\{u\})$ and $X_m(\{v\})$ can be perfectly matched.
  If $m$ is even, then after matching $X_m(\emptyset)$, 
  one configuration with $m$ particles remain, and there is
  a perfect matching for $X_m(\emptyset)$ otherwise.
  For $X_m(\{u,v\})$, we have one configuration with $m+1$ particles
  left if $m$ is odd and a perfect matching otherwise.
\qed
\end{proof}

\begin{proposition2} \label{prop:hexcyl1}
  Let $m \geq 1$. Then,
  \[
  I(C^H_{m,1}) \simeq \begin{cases}
    S^{2m/3+1} & \text{if $m \equiv_3 0$} \\
    \bullet & \text{if $m \equiv_3 1$} \\
    S^{2(m+1)/3} & \text{if $m \equiv_3 2$} \\
  \end{cases}
  \]
\end{proposition2}

\begin{figure}[htbp]
  \begin{center}
    \input{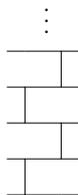}
    \caption{The graph $C^H_{4,1}$ is isomorphic to $P_9$.}
  \label{fig:proof1}
  \end{center}
\end{figure}

\begin{proof}
  It can be seen in Figure~\ref{fig:proof1}
  that the graph $C^H_{m,1}$ is isomorphic to 
  $P_{2n+1}$, the path on $2n+2$ vertices.
  The result immediately follows.
\qed
\end{proof}

\begin{proposition2} \label{prop:hexcyl2}
  Let $m \geq 1$.
  Then, the simplicial complex $I(C^H_{m,2})$ is homotopy equivalent to a 
  wedge of $F_{m+2}$ $m$-dimensional spheres.
\end{proposition2}

\begin{figure}[htbp]
  \begin{center}
    \input{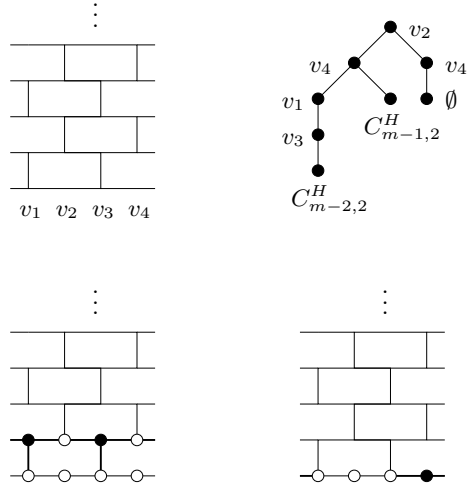}
    \caption{A recursion for $I(C^H_{m,2})$.}
  \label{fig:proof2}
  \end{center}
\end{figure}

\begin{proof}
  The Morse matching for $I(C^H_{m,2})$ is given by the matching tree in
  Figure~\ref{fig:proof2}.
  We start by splitting on vertex $v_2$.
  If there is a particle on $v_2$, then we can find a perfect matching
  on all remaining configurations by matching on $v_4$.
  If $v_2$ is empty, we split on vertex $v_4$.
  The remaining graph in the right subtree is shown in
  the lower right corner of the figure.
  We see that it is isomorphic to $C^H_{m-1,2}$,
  and that we have added one particle on the way to this node.
  In the left subtree, we eventually reach a node where the
  remaining graph is the one in the lower left corner.
  This graph is isomorphic to $C^H_{m-2,2}$ and
  the two match nodes on the way to this node adds two particles.

  The base cases when $m = 1$ and $m = 2$ are constructed by hand.
  We can construct a matching on $I(C^H_{1,2})$ so that two
  configurations remain, both with two particles.
  On $I(C^H_{2,2})$ we can similarly construct a matching
  which leaves three configurations, each with three particles.
  Using these matchings and Proposition~\ref{prop:tree}, 
  we finally construct an acyclic matching on
  $I(C^H_{m,2})$ which leaves $F_{m+2}$ unmatched configurations
  each with $m+1$ particles.
\qed
\end{proof}

As we will see in the next section when we determine some homology
groups using computer calculations, the complex $I(C^H_{5,4})$ has
non-zero homology in two dimensions.
Therefore, we can not in general hope to apply
Corollary~\ref{cor:wedgeofspheres} to the complexes $I(C^H_{m,4})$.
However, Table~\ref{tab:hexcyltrans} suggests that the
alternating sum may be simple.
We are able to find a matching tree which shows that this alternating
sum is identical to that of $I(C^H_{m,2})$ which was determined in
Proposition~\ref{prop:hexcyl2}.
This results is given in the following proposition.

\begin{proposition2} \label{prop:hexcyl4}
  Let $m \geq 1$.
  Then, $Z(C^H_{m,4}) = (-1)^{m+1} F_{m+2}.$
\end{proposition2}

\begin{figure}[htbp]
  \begin{center}
    \input{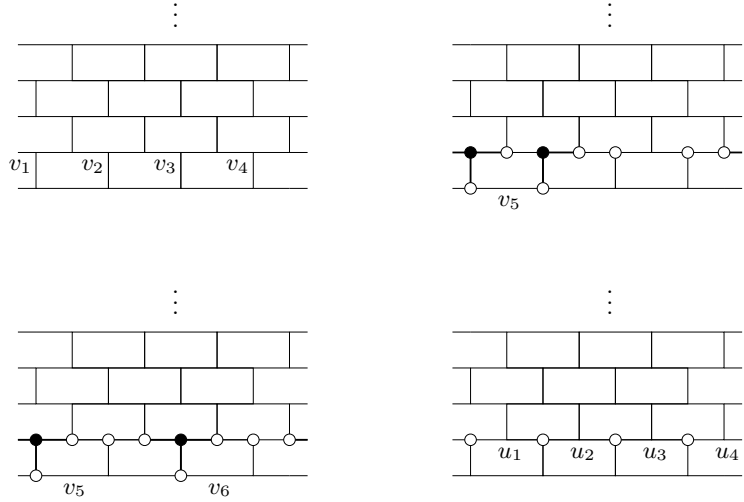}
    \caption{Split and match nodes in a matching tree for $I(C_{m,4})$.}
  \label{fig:proof3}
  \end{center}
\end{figure}

\begin{proof}
  We give a sketch of the matching tree construction used to prove
  Proposition~\ref{prop:hexcyl4}.
  We start by splitting on the vertices $v_1, \ldots, v_4$ in
  Figure~\ref{fig:proof3}.
  The first thing to note is that if particles are present
  on two consecutive 
  $v_i$ and $v_{i+1}$ (modulo 4), then we can match everything
  in this subtree. The upper right configuration shows
  this case, when for all faces $\sigma$ in this subtree, 
  $v_1, v_2 \in \sigma$.
  Here, we match on $v_5$.
  A similar situation occurs when only one $v_i$ is present.

  Therefore, we have either the case 
  $v_1, v_3 \in \sigma, v_2, v_4 \not\in \sigma$
  (or the complement), or no $v_i$ is in $\sigma$.
  The lower left configuration shows the first case, 
  where we proceed by matching on $v_5$ and $v_6$.
  In this case, which occurs twice, we
  end up with the graph $C^H_{m-2,4}$ and 4 particles
  fixed in each configuration.

  The lower right configuration shows the case when
  none of the $v_i$ are in $\sigma$.
  We then match on $u_1, \ldots, u_4$.
  In the second case, we reach the graph $C^H_{m-3,4}$
  after fixing 7 particles.
  We see here that the critical cells will be spread
  out over several dimensions.
  We therefore focus on the alternating sum, disregarding
  the homotopy type in this case.
  From the previous argument, we have the recursion
  \[
  Z(C^H_{m,4}) = 2Z(C^H_{m-2,4}) - Z(C^H_{m-3,4}),
  \]
  and base cases are
  \[
  Z(C^H_{0,4}) = -1, Z(C^H_{1,4}) = 2, Z(C^H_{2,4}) = -3,
  \]
  where $C^H_{0,4}$ is isomorphic to an 8-cycle.
  The proposition follows by solving this recursion.
\qed
\end{proof}

\subsection{The matching $M_O$}

In the rest of this section we will present the homology groups
of some independence complexes which have been computer generated.
To reduce the size of the complexes and the time of the calculations,
we first apply a particular Morse matching to each complex.
We start by introducing this matching, and the manipulations made on the
complex, before giving our results in Section~\ref{subsec:sq} and
Section~\ref{subsec:hex}

\smallskip
Let $G$ be a simple, undirected graph.
Let $O \in I(G)$ be an independent set in $G$, which we will call the
\emph{odd vertices}.
We will assume a fixed, but arbitrary total order $<$ on $O$.

The matching $M_O$ on the Hasse diagram of $P(I(G))$ is defined as follows.
The edge $(\sigma, \sigma')$ is in $M_O$ if $\sigma' \setminus \sigma = \{o\}$, where $o \in O$ is the smallest element (with respect to $<$) in $O$ for which $\sigma \cup \{o\}$ is an independent set.
According to Proposition~\ref{prop:tree}, 
the matching $M_O$ is acyclic, and therefore a Morse matching.
According to Theorem~\ref{thm:morse}, we then know that $I(G)$
is homotopy equivalent to some cell-complex,
with cells in the dimensions given by the sizes of the 
unmatched faces in $I(G)$.
We call this set of unmatched faces $X$ and note that we can
express it as
\begin{equation} \label{eq:X}
  X = \{ \sigma \in I(G) \;|\; \forall u \in O : \sigma \cap Nbd(u) \neq \emptyset \},
\end{equation}
where $Nbd(u)$ denotes the neighbourhood of $u$ in $G$.
Of course, in general, $X$ does not form a simplicial complex.
From $(\ref{eq:X})$ we see that the total order $<$ does not affect $X$.
Now, given a set $O$, and thereby a set $X$,
we define two simplicial complexes $\Delta_0$ and
$\Gamma_0$.
\[
\Delta_O = \{ \sigma \;|\; \sigma \subseteq \sigma' \in X \}, \qquad
\Gamma_O = \Delta_O \setminus X
\]

We will now prove a theorem that relates the homotopy type of the
original independence complex to that of $\Gamma_O$.
First, however, we will give a topological lemma.
In \cite{Linusson06} it was stated for the case when $\Delta$ is a simplex.
This slightly more general version is needed for the proof of our theorem.

\begin{lemma2} \label{lem:suspension}
  Let $\Gamma \subseteq \Delta$ be simplicial complexes and assume that 
  $\Delta$ is contractible.
  Then, $\Delta/\Gamma \simeq \mathrm{susp}(\Gamma)$.
\end{lemma2}

\begin{proof}
  The concept of a \emph{mapping cone} is described in 
  \cite[Chapter~0]{Hatcher02}.
  In our case, the mapping cone is the simplicial complex
  $X = \Delta \sqcup_{\mathrm{id}} \mathrm{cone}(\Gamma$),
  i.e, we raise a cone in $\Delta$ over the subcomplex $\Gamma$.
  Note that $X/\mathrm{cone}(\Gamma) \simeq \Delta/\Gamma$ and that
  $X/\Delta \simeq \mathrm{susp}(\Gamma)$.
  Now, since both $\mathrm{cone}(\Gamma)$ and $\Delta$ are contractible,
  we have $X/\mathrm{cone}(\Gamma) \simeq X \simeq X/\Delta$ which
  finishes the proof.
\qed
\end{proof}

We are now ready to state and prove the following theorem, which allows us
to simplify our calculations of the homology groups.

\begin{theorem} \label{thm:matching}
  Let $O \in I(G)$ be an independent set and assume that $\Delta_O$ is contractible.
  Then,
  \begin{equation}
    I(G) \simeq \mathrm{susp}(\Gamma_O).
  \end{equation}
\end{theorem}

\begin{proof}
  The set $X$ is given by $(\ref{eq:X})$.
  Another way to express this is to say that $\sigma \in X$ if and only if
  $O \subseteq \bigcup_{u \in \sigma} Nbd(u)$,
  From the latter definition we can immediately see that
  if $\sigma \subseteq \sigma'$ and $\sigma \in X$, then we have $\sigma' \in X$.
  Let $\Sigma_O = I(G) \setminus X$.
  The previous argument shows that $\Sigma_O$ is a simplicial complex.
  Figure~\ref{fig:mo} illustrates the relationship between $X$ and the simplicial complexes $\Delta_O, \Gamma_O$ and $\Sigma_O$.
  \begin{figure}[htbp]
    \begin{center}
      \input{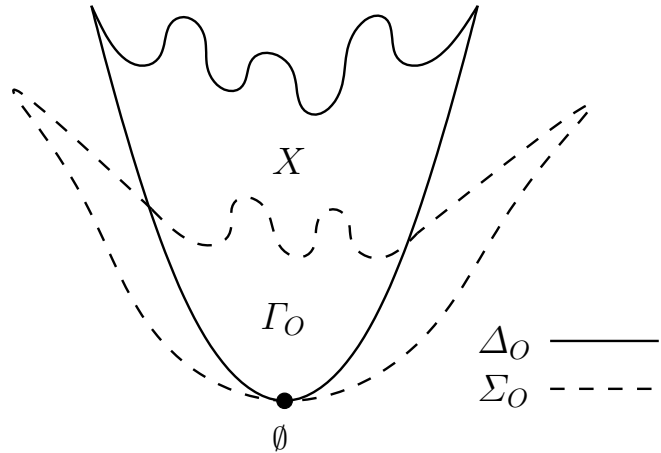}
      \caption{The set $X$ and the simplicial complexes $\Delta_O, \Gamma_O$ and $\Sigma_O$}
      \label{fig:mo}
    \end{center}
  \end{figure}

  Since $M_O$ gives a perfect matching on $\Sigma_O$, we can collapse
  $\Sigma_O$ to a point in $I(G)$.
  But this is the same as identifying $\Gamma_O$ with a point in $\Delta_O$.
  Thus,
  \[
  I(G) \simeq I(G)/\Sigma_O \simeq \Delta_O/\Gamma_O.
  \]
  Using Lemma~\ref{lem:suspension}, we now see that $I(G) \simeq \mathrm{susp}(\Gamma_0)$. 
\qed
\end{proof}

The reason for using $\Gamma_O$ instead of $I(G)$ is two-fold.
First, we have used a matching to reduce the size of the complex, and
thereby, hopefully, reduced the time to compute its homology.
Secondly, in many cases a description of the complex $\Gamma_O$ is easier
to obtain and more compact than a description of the original complex.

\subsection{Graphs on the square grid} \label{subsec:sq}

Let $G = C_{m,n}$ and assume that $n$ is even.
In this case we choose the independent set $O = \{ u = (i,j) \in V(G) \;|\; \text{$i+j$ is odd} \}$.
Then, $\Delta_O \cong \Delta_{mn/2-1}$ and thus contractible.
The complex $\Gamma_O$ can be generated from the set $O$ in the following way.
For each $u \in O$ we have a maximal face $\{ v \in V(G) \setminus O \;|\; v \not\in Nbd(u) \} \in \Gamma_O$.
We see that for the maximal faces $\sigma \in \Gamma_O$, either $\dim (\sigma) = m n / 2 - 5$ or $\dim (\sigma) = m n /2 - 4$, depending on whether the corresponding $v \in V(G) \setminus O$ is on the (non-identified) boundary of the cylinder or not.

\begin{example}
  Let $m = 3$, $n = 4$, and
  $O = \{ (1,2), (1,4), (2,1), (2,3), (3,2), (3,4) \}$.
  Then, the six maximal faces of $\Gamma_O$ are as shown in
  Figure~\ref{fig:ex1}
  and the complex itself can be realised as in Figure~\ref{fig:ex2}.
  We can see that it is homotopy equivalent to a wedge of three 
  $1$-dimensional spheres.
  Therefore, we conclude from Theorem~\ref{thm:matching} that 
  $I(C_{3,4}) \simeq \bigvee_3 S^2$.
  \begin{figure}[htbp]
    \begin{center}
      \input{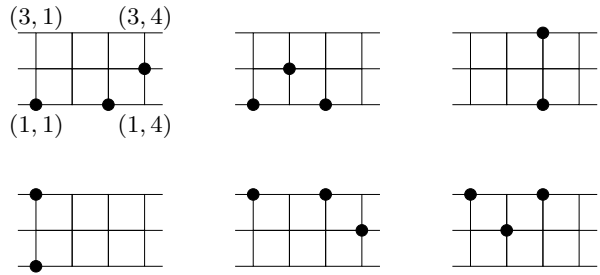}
      \caption{The six maximal faces of $\Gamma_O$.}
      \label{fig:ex1}
    \end{center}
  \end{figure}
  \begin{figure}[htbp]
    \begin{center}
      \input{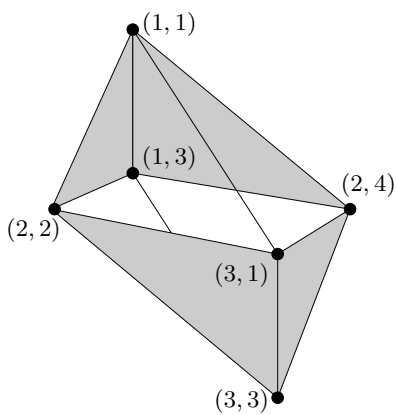}
      \caption{The simplicial complex $\Gamma_{O} \simeq \bigvee_{3} S^{1}$.}
      \label{fig:ex2}
    \end{center}
  \end{figure}
\end{example}

For odd $n$, we let $O = \{ u = (i,j) \in V(G) \;|\; j \neq n, \text{$i+j$ is odd} \}$.
The complex $\Delta_O$ is still contractible but the
maximal faces of $\Gamma_O$ are less regular.
\medskip

Table~\ref{tab:sqcyleven} and Table~\ref{tab:sqcylodd} show 
the calculated homology groups for the complex $I(C_{m,n})$ 
for small $m$ and $n$.
An entry $(k,d)$ in the table should be read as $H_k(I(C_{m,n})) \cong \mathbb{Z}^d$.
Missing entries are cases where Polymake ran for too long, 
or more often, ran out of memory.
The homology vanishes in all dimensions that are not listed.
Note that for $n \leq 5$, the homology groups in 
these tables can be derived, for arbitrary $m$,
 from the results in Section~\ref{subsec:small}.
For odd $n$, we observe an interesting symmetry, which we pose as a conjecture.

\begin{conjecture}
  For all $i \geq 0$ and $j, k \geq 1$, the following holds.
  \[
  H_i(I(C_{j,2k+1})) \cong H_i(I(C_{k,2j+1}))
  \]
\end{conjecture}

Table~\ref{tab:sqcyltrans} is taken from \cite{Jonsson06}
and shows the values of $Z(C_{m,n})$ for small $m$ and $n$,
calculated using transfer matrices.
The cases when $n$ is odd are left out as they all satisfy
Conjecture~\ref{conj:cylinder}.
In these cases, it was also observed, and conjectured to hold in
general, that all the roots of the characteristic polynomials
were roots of unity.
In particular, this implies that the rows of 
Table~\ref{tab:sqcyltrans} are periodic.

\begin{table}[htbp]
  \begin{center}
    \begin{tabular}{|c|r||c|c|c|c|c|c|c|c|}
      \hline
      \multicolumn{2}{|c||}{\multirow{2}{*}{}} & \multicolumn{8}{|c|} m  \\
      \cline{3-10}
      \multicolumn{2}{|c||}{}&1&2&3&4&5&6&7&8 \\
      \hline
      \hline
      \multirow{7}{*}{\, n \,}
      & 2 & (0,1) & (0,1) & (1,1) & (1,1) & (2,1) & (2,1) & (3,1) & (3,1) \\
      \cline{2-10}
      & 4 & (0,1) & (1,3) & (2,3) & (3,5) & (4,5) & (5,7) & (6,7) & (7,9) \\
      \cline{2-10}
      & 6 & (1,2) & (2,1) & (3,1) & (5,4) & (6,1) & (7,1), (8,2) &  (9,4) & (11,7) \\
      \cline{2-10}
      & 8 & (2,1) & (3,3) & (5,5) & (7,5) & (8,1), (9,4) & (11,7) & \multicolumn{2}{r|}{} \\
      \cline{2-8}
      & 10 & (2,1) & (4,1) & (7,1) & (8,1), (9,2) & \multicolumn{4}{r|}{} \\
      \cline{2-6}
      & 12 & (3,2) & (5,3) & (8,3) & (11, 8) & \multicolumn{4}{r|}{} \\
      \cline{2-6}
      & 14 & (4,1) & (6,1) & (9,1) & \multicolumn{5}{r|}{} \\
      \hline
    \end{tabular}
  \end{center}
  \caption{Homology of $I(C_{m,n})$ for even $n$.}
  \label{tab:sqcyleven}
\end{table}

\begin{table}[htbp]
  \begin{center}
    \begin{tabular}{|c|r||c|c|c|c|c|c|c|c|c|c|c|}
      \hline
      \multicolumn{2}{|c||}{\multirow{2}{*}{}} & \multicolumn{11}{|c|} m \\
      \cline{3-13}
      \multicolumn{2}{|c||}{}&1&2&3&4&5&6&7&8&9&10&11 \\
      \hline
      \hline
      \multirow{9}{*}{\, n \,}
      &3 & (0,2) & (1,1) & (1,1) & (2,2) & (3,1) & (3,1) & (4,2) & (5,1) & (5,1) & (6,2) & (7,1) \\
      \cline{2-13}
      &5 & (1,1) & (1,1) & (3,1) & (3,1) & (5,1) & (5,1) & (7,1) & (7,1) & (9,1) & (9,1) & \\
      \cline{2-12}
      &7 & (1,1) & (3,1) & (5,1) & (5,1) & (7,1) & (9,1) & (11,1) & \multicolumn{4}{r|}{} \\
      \cline{2-9}
      &9 & (2,2) & (3,1) & (5,1) & (7,1), (8,3) & (9,1) & \multicolumn{6}{r|}{} \\
      \cline{2-7}
      &11 & (3,1) & (5,1) & (7,1) & (9,1) & \multicolumn{7}{r|}{} \\
      \cline{2-6}
      &13 & (3,1) & (5,1) & (9,1) & \multicolumn{8}{r|}{} \\
      \cline{2-5}
      &15 & (4,2) & (7,1) & (11,1) & \multicolumn{8}{r|}{} \\
      \cline{2-5}
      &17 & (5,1) & (7,1) & \multicolumn{9}{r|}{} \\
      \cline{2-4}
      &19 & (5,1) & (9,1) & \multicolumn{9}{r|}{} \\
      \hline
    \end{tabular}
  \end{center}
  \caption{Homology of $I(C_{m,n})$ for odd $n$.}
  \label{tab:sqcylodd}
\end{table}

\begin{table}[htbp]
  \begin{center}
    \begin{tabular}{|c|r||r|r|r|r|r|r|r|r|r|r|r|}
      \hline
      \multicolumn{2}{|c||}{\multirow{2}{*}{$Z(C_{m,n})$}} & \multicolumn{11}{|c|} m \\
      \cline{3-13}
      \multicolumn{2}{|c||}{}&\makebox[16pt][r]{1}&\makebox[16pt][r]{2}&\makebox[16pt][r]{3}&\makebox[16pt][r]{4}&\makebox[16pt][r]{5}&\makebox[16pt][r]{6}&\makebox[16pt][r]{7}&\makebox[16pt][r]{8}&\makebox[16pt][r]{9}&\makebox[16pt][r]{10}&\makebox[16pt][r]{11} \\
      \hline
      \hline
      \multirow{11}{*}{\,\, n \,\,}
        & 2 & -1 & -1 &  1 &  1 & -1 & -1 &  1 &  1 & -1 & -1 &  1 \\
      \cline{2-13}
        & 4 & -1 &  3 & -3 &  5 & -5 &  7 & -7 &  9 & -9 & 11 &-11 \\
      \cline{2-13}
        & 6 &  2 & -1 &  1 &  4 & -1 & -1 &  4 &  1 & -1 &  2 &  1 \\
      \cline{2-13}
        & 8 & -1 &  3 &  5 &  5 &  3 &  7 &  1 &  1 & -1 &  3 & -3 \\
      \cline{2-13}
        &10& -1 & -1 &  1 &  1 &  9 & -1 &  1 &  1 &-11 & -1 &  1 \\
      \cline{2-13}
        &12&  2 &  3 & -3 &  8 & -5 &  7 &  8 &  9 & -9 & 14 &-11 \\
      \cline{2-13}
        &14& -1 & -1 &  1 &  1 & -1 & 13 &  1 &  1 & 13 & -1 & 15 \\
      \cline{2-13}
        &16& -1 &  3 &  5 &  5 &  3 &  7 &  1 & 33 & -1 &  3 & 13 \\
      \cline{2-13}
        &18&  2 & -1 &  1 &  4 & -1 & -1 & 22 &  1 & -1 & 38 &  1 \\
      \cline{2-13}
        &20& -1 &  3 & -3 &  5 &  5 &  7 & -7 &  9 & 41 & 11 &-11 \\
      \cline{2-13}
        &22& -1 & -1 &  1 &  1 & -1 & -1 &  1 & 23 & -1 & -1 & 89 \\
      \cline{2-13}
        &24&  2 &  3 &  5 &  8 &  3 &  7 & 16 &  1 & -1 & 78 & -3 \\
      \hline
    \end{tabular}
  \end{center}
  \caption{$Z(C_{m,n})$ for small $m$ and $n$, when $n$ is even. 
    Taken from Jonsson \cite{Jonsson06}.}
  \label{tab:sqcyltrans}
\end{table}

\smallskip
Both for even and odd $n$, we find cases where there is non-zero homology
in more than one dimension.
This means that in these cases there is no easy way to determine the 
alternating sums via Morse matchings and the homotopy type.
In the odd $n$-case, the \emph{bad} example shows up for the graph
$C_{4,9}$, where $Z(C_{4,9}) = -2$.
It is therefore still possible that the independence complexes in the
other odd $n$-cases, where $Z(C_{m,n})$ is conjectured to be $1$, are simple.

\subsection{Graphs on the hexagonal grid} \label{subsec:hex}

Take any finite subgraph $H'$ of the hexagonal grid $H$.
As in the case of the square grid with even $n$, there is a natural
notion of \emph{even} and \emph{odd} vertices given by the parity
of the vertices.
Let $G$ be the graph $H'$, possibly with some identifications made
on vertices of the same parity
and let $O$ be the set of odd vertices in $G$,
where the parity is inherited from $H$.
Then, $\Delta_O$ is a simplex of dimension $|V(G)|-|O(G)|-1$ and thus contractible.
The complex $\Gamma_O$ is determined in the same way as was done for the
square grid with even $n$.
For each $u \in O$ we have a maximal face $\{ v \in V(G) \setminus O \;|\; v \not\in Nbd(u) \} \in \Gamma_O$.
\smallskip

Table~\ref{tab:hexcyl} shows the homology of the complex $I(C^H_{m,n})$ for small $m$ and $n$. As before, an entry $(k,d)$ means that $H_k(I(C^H_{m,n})) \cong \mathbb{Z}^d$ and that the homology vanishes in all other dimensions.
As noted in the square grid case, the entries for $n=2$ could be derived,
for arbitrary $m$ from Proposition~\ref{prop:hexcyl2}.
We note with interest that there seems to be a relation between the
dimensions in which homology appears for 
$C^H_{i,j}$ and $C^H_{j-1,i+1}$, but the data is quite limited.
Table~\ref{tab:hexcyltrans} shows some of the first values of the alternating sum $Z(C^H_{m,n})$. 
They have been generated using transfer matrices.
We have not found any nice factorisation of these polynomials and
they seem to have mostly roots which are not roots of unity.

\begin{table}[htbp]
  \begin{center}
    \begin{tabular}{|c|r||c|c|c|c|c|c|c|}
      \hline
      \multicolumn{2}{|c||}{\multirow{2}{*}{}} & \multicolumn{7}{|c|} m \\
      \cline{3-9}
      \multicolumn{2}{|c||}{}&2&3&4&5&6&7&8\\
      \hline
      \hline
      \multirow{7}{*}{\, n \,}
      &2 & (2,3) & (3,5) & (4,8) & (5,13) & (6,21) & (7,34) & (8,55) \\
      \cline{3-9}
      &3 & (4,5) & (6,7) & (7,3) & (9,22) & (11,23) & (12,24) &\\
      \cline{3-8}
      &4 & (6,3) & (7,4), (9,1) & (10,8) & (11,8), (13,5) & \multicolumn{3}{r|}{} \\
      \cline{3-6}
      &5 & (7,6) & (10,11) &  \multicolumn{5}{r|}{} \\  
      \cline{3-4}
      &6 & (9,15) & (11,4), (13,13) &  \multicolumn{5}{r|}{} \\  
      \cline{3-4}
      &7 & (11,8) &  \multicolumn{6}{r|}{} \\  
      \cline{3-3}
      &8 & (12,19) & \multicolumn{6}{r|}{} \\ \hline
    \end{tabular}
  \end{center}
  \caption{Homology of $I(C^H_{m,n})$.}
  \label{tab:hexcyl}
\end{table}

\begin{table}[htbp]
  \begin{center}
    \begin{tabular}{|c|r||r|r|r|r|r|r|r|r|r|r|r|r|}
      \hline
      \multicolumn{2}{|c||}{\multirow{2}{*}{$Z(C^H_{m,n})$}} & \multicolumn{12}{|c|} m \\
      \cline{3-14}
      \multicolumn{2}{|c||}{}&\quad1&\quad2&\quad3&\quad4&\quad5&\quad6&\quad7&\quad8&\quad9&\quad10&\quad11&\quad12\\
      \hline
      \hline
      \multirow{9}{*}{\,\, n \,\,}
      &1 & 0 & 1 & -1 & 0 & 1 & -1 & 0 & 1 & -1 & 0 & 1 & -1 \\
      \cline{3-14}
      &2 & 2 & -3 & 5 & -8 & 13 & -21 & 34 & -55 & 89 & -144 & 233 & -377 \\
      \cline{3-14}
      &3 & 0 & -5 & -7 & 3 & 22 & 23 & -24 & -92 & -67 & 141 & 367 & 152 \\
      \cline{3-14}
      &4 & 2 & -3 & 5 & -8 & 13 & -21 & 34 & -55 & 89 & -144 & 233 & -377 \\
      \cline{3-14}
      &5 & 0 & 6 & -11 & -5 & 51 & -76 & -60 & 416 & -536 & -655 & 3351 & -3646 \\
      \cline{3-14}
      &6 & 2 & 15 & 17 & 55 & 160 & 231 & 886 & 1664 & 3947 & 11121 & 21065 & 59296 \\
      \cline{3-14}
      &7 & 0 & 8 & -15 & -35 & 57 & 34 & -42 & 687 & 20 & -4207 & -2379 & 3611 \\
      \cline{3-14}
      &8 & 2 & -19 & 37 & -88 & 533 & -725 & 3466 & -11927 & 21417 & -105552 & 273881 & -682665 \\
      \cline{3-14}
      &9 & 0 & -41 & -43 & 183 & 958 & 941 & -9924 & -22943 & 19265 & 289806 & 587437 & -1949599 \\
      \hline
    \end{tabular}
  \end{center}
  \caption{$Z(C^H_{m,n})$ for small $m$ and $n$.}
  \label{tab:hexcyltrans}
\end{table}

In the case of the hexagonal lattice, we have also calculated the
homology for the toroidal case, i.e, for the complex $I(T^H_{m,n})$.
Table~\ref{tab:hextorus} shows the homology for small cases.
Table~\ref{tab:hextorusmatrices} lists $Z(C^T_{m,n})$ for some
small values of $m$ and $n$ calculated using the transfer matrices.
For $m \leq 4$, the characteristic polynomials $P(m)$ of these matrices
factor reasonably well
After removing powers of $t$, the polynomials are given by
\[
\begin{tabular}{rcl}
  $P(1)$ & = & $1 + t + t^2$ \\
  $P(2)$ & = & $1+t-t^2$ \\
  $P(3)$ & = & $(1-t+t^2)(1-t+2t^2+t^3)$ \\
  $P(4)$ & = & $(1+t-t^2)(1-t)^3(1+t)^3(1+t^2)$ \\
\end{tabular}
\]
For larger values of $m$, the polynomials do not factor as nicely,
and seem to have mostly roots which are not roots of unity.

\begin{table}[htbp]
  \begin{center}
\scalebox{0.95}{
    \begin{tabular}{|c|r||c|c|c|c|c|c|c|}
      \hline
      \multicolumn{2}{|c||}{\multirow{2}{*}{}} & \multicolumn{7}{|c|} m \\
      \cline{3-9}
      \multicolumn{2}{|c||}{}&2&3&4&5&6&7&8\\
      \hline
      \hline
      \multirow{3}{*}{\, n \,}
      &2 & (1,3) & (2,4) & (3,7) & (4,11) & (5,18) & (6,29) & (7,47) \\
      \cline{3-9}
      &3 & (2,4) & (4,10) & (6,4) & (7,17) & (9,32) & (10,1), (11,3) & 
      (12,76) \\
      \cline{3-9}
      &4 & (3,7) & (6,4) & (7,15) & (9,1), (10,12) & (11,20), (12,1),
      (13,3) & \multicolumn{2}{r|}{} \\
      \hline
    \end{tabular}
}
  \end{center}
  \caption{Homology of $I(T^H_{m,n})$.}
  \label{tab:hextorus}
\end{table}

\begin{table}[htbp]
  \begin{center}
    \begin{tabular}{|c|r||r|r|r|r|r|r|r|}
      \hline
      \multicolumn{2}{|c||}{\multirow{2}{*}{$Z(T^H_{m,n})$}} & \multicolumn{7}{|c|} m \\
      \cline{3-9}
      \multicolumn{2}{|c||}{}&\quad1&\quad2&\quad3&\quad4&\quad5&\quad6&\quad7\\
      \hline
      \hline
      \multirow{8}{*}{\,\, n \,\,}
      &1 & -1 & -1 & 2 & -1 & -1 & 2 & -1 \\
      \cline{3-9}
      &2 & -1 & 3 & -4 & 7 & -11 & 18 & -29 \\
      \cline{3-9}
      &3 & 2 & -4 & -10 & -4 & 17 & 32 & 2 \\
      \cline{3-9}
      &4 & -1 & 7 & -4 & 15 & -11 & 22 & -29 \\
      \cline{3-9}
      &5 & -1 & -11 & 17 & -11 & -51 & 127 & -36 \\
      \cline{3-9}
      &6 & 2 & 18 & 32 & 22 & 127 & 192 & 394 \\
      \cline{3-9}
      &7 & -1 & -29 & 2 & -29 & -36 & 394 & 552 \\
      \cline{3-9}
      &8 & -1 & 47 & -76 & 55 & -411 & 1478 & 83 \\
      \hline
    \end{tabular}
  \end{center}
  \caption{$Z(T^H_{m,n})$ for small $m$ and $n$.}
  \label{tab:hextorusmatrices}
\end{table}




\section{A second matching} \label{sec:matching}

In this section we introduce an alternative matching on the face poset of the complex $I(C_{m,n})$.
While in Section~\ref{sec:explicit}, the matching presented was always acyclic, we will make no such claims here.
This means that we can not rely on the discrete Morse theory for our results.
Still, we will be able to derive results on the alternating sum $Z(C_{m,n})$.
A similar idea was used with great success on the hard square model on a torus in \cite{Jonsson06}.

\medskip
Let $M$ be any matching on the underlying graph of $P(I(C_{m,n}))$.
Then, for all $(\sigma, \tau) \in M$, we have $(-1)^{|\sigma|} + (-1)^{|\tau|} = 0$ and therefore, 
\[
\sum_{\sigma \in I(C_{m,n})} (-1)^{|\sigma|} = \sum_{\sigma \in X} (-1)^{|\sigma|},
\]
where we let $X$ denote the set of unmatched faces in $I(C_{m,n})$.
Thus, we can calculate $Z(C_{m,n})$ using only
the faces in the set $X$.
\medskip

In order to describe our matching we need to introduce
some notation and terminology.
Let $\sigma$ be a face in $I(C_{m,n})$.
Define
\[
\pi(\sigma) := \{ (1,j) \in \sigma \}.
\]
Let $\pi(\sigma) = \{x_1, \ldots, x_k\}$ be ordered by the
second component of the $x_i$.
For $1 \leq i \leq k$, we define the \emph{interval} $N_i$
of $\sigma$ at $x_i = (1,j_i)$
to be the set
\[
N_i = \{(1,j_{i-1}+1), (1,j_{i-1}+2), \ldots, (1,j_i) \},
\]
where the addition is modulo $n$ in the second component,
and the indices $0$ and $k$ are identified.
Let $N(\sigma) = \{ N_1, \ldots, N_k \}$ be the set of all
intervals.
We call the sequence $[|N_1|, \ldots, |N_k|]$ the \emph{signature}
of $\sigma$ and we identify signatures which are identical up to
cyclic rotations.
We distinguish between \emph{even} and \emph{odd intervals} by the parity of
their sizes and between \emph{even} and \emph{odd positions} in an interval
$N_i$ by saying that the element $x_i$ is in an even (odd) position
if $N_i$ is an even (odd) interval,
and extending this notion of parity to the rest of the interval.
Let $\pi(\sigma) = \pi_o(\sigma) \cup \pi_e(\sigma)$, where $\pi_o$ is the set of particles in odd position, and $\pi_e(\sigma)$ the set of particles in even position.

  \begin{figure}[htbp]
    \begin{center}
      \input{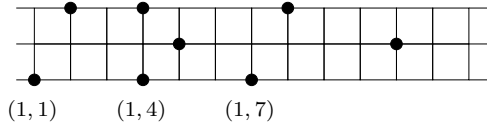}
      \caption{A face $\sigma \in C_{3,13}$.}
      \label{fig:interval}
    \end{center}
  \end{figure}

\begin{example} \label{ex:1}
  Let $\sigma \in C_{3,13}$
  be the configuration in Figure~\ref{fig:interval}.
  Here $\pi(\sigma) = \pi_o(\sigma) = \{(1,1),(1,4),(1,7)\}$ and the intervals $N(\sigma) = \{N_1, N_2, N_3\}$ are given by
  \[
  \begin{tabular}{rcl}
    $N_1$& = &$\{ (1,8), (1,9), (1,10), (1,11), (1,12), (1,13), (1,1) \}$, \\
    $N_2$& = &$\{ (1,2), (1,3), (1,4) \}$, and \\
    $N_3$& = &$\{ (1,5), (1,6), (1,7) \}$. \\
  \end{tabular}
  \]
  In the interval $N_3$, the positions $(1,5)$ and $(1,7)$ are odd, while in $N_1$, the odd positions are $(1,8)$, $(1,10)$, $(1,12)$ and $(1,1)$.
  The signature of $\sigma$ is $[7,3,3]$ which is the same as $[3,7,3]$ or $[3,3,7]$.
\end{example}

We start by defining a closure operator $\sigma \mapsto {\widehat \sigma}$
on the poset $P(I(C_{m,n}))$.
The idea is to use the equivalence relation
defined by $\sigma \sim \tau$ if and only if
${\widehat \sigma} = {\widehat \tau}$ to group faces
in equivalence classes.
By constructing matchings on the non-singleton classes of $\sim$
we can then restrict our study to a particular 
collection of singleton classes.
As before, let $\pi(\sigma) = \{x_1, \ldots, x_k\}$.
The face ${\widehat \sigma}$ is defined as follows.
\begin{itemize}
\item $\sigma \subseteq {\widehat \sigma}$.
\item For each $i = 1,\ldots,k$, if $x_i$ is in an even position and if $x_i$ is free in $\sigma$, then $x_i \in {\widehat \sigma}$.
\end{itemize}

Let $\sim$ be the equivalence relation described above and let
$X_\sigma$ be the equivalence class of $\sigma$ with respect to $\sim$.
That is, 
\[
X_\sigma = \{ \tau \in I(C_{m,n}) \;|\; \tau \sim \sigma \}.
\]
Let $\sigma$ be the face from Example~\ref{ex:1}. Then,
${\widehat \sigma} = \sigma + (1,9)$
and $X_\sigma = \{\sigma, \sigma + (1,9\}$.

\begin{lemma2} \label{lemma:structure}
  Let $\sigma \in I(C_{m,n})$.
  \begin{itemize}
  \item Let $x$ be a free and even position in an interval $N_i \in N(\sigma)$.
    Then, $\sigma \sim \sigma+x$, $\pi_o(\sigma+x) = \pi_o(\sigma)$, and $\pi_e(\sigma+x) = \pi_e(\sigma) + x$.
  \item Assume that $|N(\sigma)| > 1$. Let $x_i \in \pi_e(\sigma)$.
    Then, $\sigma \sim \sigma-x_i$, $\pi_o(\sigma-x_i) = \pi_o(\sigma)$, and $\pi_e(\sigma-x_i) = \pi_e(\sigma) - x_i$.
  \end{itemize}
\end{lemma2}

\begin{proof}
  The first part follows from the observation that adding the element $x$ to
  $\sigma$ divides the interval $N_i$ into one even interval, with last
  position in $x$, and one interval of the same parity as $N_i$, given by
  the position $x_i$.

  The second part holds since gluing an even interval with its following
  interval (which is necessarily different due to the lower
  bound condition on $N(\sigma)$) 
  creates an interval of the same parity as the latter one.
  The parity of the positions in the new interval are the same as in
  the old intervals
  and the vacated position $x_i$ is free and even in $\sigma-x_i$,
  which shows that 
  $\sigma \sim \sigma-x_i$.
\qed
\end{proof}

We partition the faces $\sigma \in I(C_{m,n})$ into three sets 
$P_1$, $P_2$ and $P_3$, depending on the signature 
and the size of $X_\sigma$.
\begin{itemize}
  \item $P_1$ is the set of $\sigma$ for which $|X_\sigma| > 1$ and for
    which there is at least one odd interval in $N(\sigma)$.
  \item $P_2$ is the set of $\sigma$ for which $|X_\sigma| > 1$ and 
    all the intervals in $N(\sigma)$ are even.
  \item $P_3$ is the set of $\sigma$ for which $X_\sigma$ is a singleton.
\end{itemize}
Note that Lemma~\ref{lemma:structure} implies that if $\sigma \in P_i$,
then $\tau \in P_i$ for all $\tau \in X_\sigma$, $i \in \{1, 2, 3\}$.

We can use Lemma~\ref{lemma:structure} to describe the sets $X_\sigma$
more explicitly.
The idea is to use the first part to move from 
$\sigma$ to ${\widehat \sigma}$.
In the process, we add only even intervals.
Then, we can use the second part to move from ${\widehat \sigma}$ to any
$\tau \in X_\sigma$ by doing the reverse operation of gluing an even
interval with its following interval.

\begin{cor}
  Assume that $\sigma \in P_1$.
  Then, $\tau \in X_\sigma$ if and only if $\tau$ can be obtained from
  ${\widehat \sigma}$ by removing a subset of $\pi_e({\widehat \sigma})$.
\end{cor}

\begin{cor}
  Assume that $\sigma \in P_2$.
  Then, $\tau \in X_\sigma$ if and only if $\tau$ can be obtained from
  ${\widehat \sigma}$ by removing a proper subset of $\pi_e({\widehat \sigma})$.
\end{cor}

We can therefore describe the structure of the equivalence classes in the
following way.
\begin{equation} \label{eq:structure}
X_\sigma = \begin{cases}
  \{{\widehat \sigma} \setminus S \;|\; S\subseteq \pi_e({\widehat \sigma})\} &
  \text{if $\sigma \in P_1$,} \\
  \{{\widehat \sigma} \setminus S \;|\; S\varsubsetneq \pi({\widehat \sigma})\} &
  \text{if $\sigma \in P_2$,} \\
  \{\sigma\} & \text{if $\sigma \in P_3$.} \\
\end{cases}
\end{equation}


Using the description of $X_\sigma$ in (\ref{eq:structure}), we can
show that the equivalence classes in $P_1$ and $P_2$ behave well.
These results are given in Lemma~\ref{lemma:perfect} and 
Lemma~\ref{lemma:almostperfect}.

\begin{lemma2} \label{lemma:perfect}
  Let $\sigma \in P_1$.
  Then,
  \[
  \sum_{\tau \in X_\sigma} (-1)^{|\tau|} = 0.
  \]
\end{lemma2}

\begin{proof}
  Since $|X_\sigma| > 1$, we know that $\pi_e({\widehat \sigma})$ is non-empty.
  In particular, we can choose an $x_i \in \pi_e({\widehat \sigma})$ and match
  $\tau - x_i$ with $\tau + x_i$ for each $\tau \in X_\sigma$.
\qed
\end{proof}

\begin{lemma2} \label{lemma:almostperfect}
  Let $\sigma \in P_2$. Then,
  \[
  \sum_{\tau \in X_\sigma} (-1)^{|\tau|} = (-1)^{|\sigma \setminus \pi(\sigma)|+1}.
  \]
\end{lemma2}

\begin{proof}
  Note that for all $\sigma \in P_2$, we have
  $\pi_e({\widehat \sigma}) = \pi({\widehat \sigma})$.
  According to (\ref{eq:structure}), 
  the set $X_\sigma$ consists of the faces ${\widehat \sigma} \setminus S$
  where $S \varsubsetneq \pi({\widehat \sigma})$.
  Therefore,
  \[
  \sum_{\tau \in X_\sigma} (-1)^{|\tau|} =
  \left( \sum_{\tau = {\widehat \sigma} \setminus S, S \subseteq \pi({\widehat \sigma})} (-1)^{|\tau|} \right) - (-1)^{|{\widehat \sigma} \setminus \pi({\widehat \sigma})|} = 0 - (-1)^{|\sigma \setminus \pi(\sigma)|}.
  \]
\qed
\end{proof}

Lemma~\ref{lemma:almostperfect} needs a bit more work to use than
Lemma~\ref{lemma:perfect}.
In order to find the sum of $(-1)^{|\sigma|}$ over all $\sigma \in P_2$,
we must use Lemma~\ref{lemma:almostperfect} on each equivalence class
in $P_2$.
Luckily, there is an easy way to describe all equivalence classes.
For $\sigma \in P_2$, let $Y_\sigma = \{ \tau \in P_2 \;|\; \tau \setminus \pi(\tau) = \sigma \setminus \pi(\sigma) \}$.
That is, $Y_\sigma$ is the set of faces in $P_2$ which are equal to $\sigma$
on rows 2 through $m$.
Then, since all faces in $P_2$ contain only even intervals,
$Y_\sigma$ can be partitioned into two sets,
\[
Y^1_\sigma = \{ \tau \in Y_\sigma \;|\; \pi(\tau) \subseteq 
\{(1,1),(1,3),\ldots,(1,n-1)\} \}
\]
and
\[
Y^2_\sigma = \{ \tau \in Y_\sigma \;|\; \pi(\tau) \subseteq 
\{(1,2),(1,4),\ldots,(1,n)\} \}
\]
Now, the set of non-empty $Y^1_\sigma$ and $Y^2_\sigma$ equals
the set of equivalence classes in $P_2$.
This observation leads to the following lemma.

\begin{lemma2} \label{lemma:p2}
  \[
  \sum_{\sigma \in P_2} (-1)^{|\sigma|} = 
  - 2 \sum_{\tau \in I(C_{m-1,n})} (-1)^{|\tau|} + 
  \begin{cases}
    2 \cdot (-1)^{mn/4} & \text{ if $m$ is even,} \\
    0 & \text{ if $m$ is odd.} \\
  \end{cases}
  \]
\end{lemma2}

\begin{proof}
  Following the previous discussion, we reorder the sum to be over
  the sets $Y_\sigma$ and since $Y_\sigma$ is uniquely determined
  by rows $2$ through $m$ of $\sigma$, we will sum over
  faces $\mu \in C_{m-1,n}$, abusing notation slightly to write
  $Y_\mu = Y_\sigma$ when $\mu$ equals $\sigma \setminus \pi(\sigma)$
  as sets.
  \[
    \sum_{\sigma \in P_2} (-1)^{|\sigma|} =
    \sum_{Y_\sigma} \sum_{\tau \in Y_\sigma} (-1)^{|\tau|} =
    \sum_{\mu \in C_{m-1,n}} \left(
    \sum_{\tau \in Y^1_{\mu}} (-1)^{|\tau|} +
    \sum_{\tau \in Y^2_{\mu}} (-1)^{|\tau|}
    \right)
  \]
  When $Y^1_{\mu}$ and $Y^2_{\mu}$ in the inner sum are
  non-empty, we can use Lemma~\ref{lemma:almostperfect}.
  The set $Y^1_{\mu}$ is empty precisely when the first row
  of $\mu$ (which is the second row of the corresponding $\sigma$)
  has particles in positions $(1,1),(1,3),\ldots,(1,n-1)$.
  Let $X_1$ be the set of such $\mu$ and let $X_2$ be the set
  of $\mu$ which has particles in positions $(1,2),(1,4),\ldots,(1,n)$.
  This gives us the following expression.
  \[
    \sum_{\sigma \in P_2} (-1)^{|\sigma|} =
    \sum_{\mu \in C_{m-1,n}} 2 \cdot (-1)^{|\mu|+1} 
    - \sum_{\tau \in X_1} (-1)^{|\tau|+1} 
    - \sum_{\tau \in X_2} (-1)^{|\tau|+1}
  \]
  The lemma follows by constructing a matching tree for the set $X_1$
  (or equivalently, for $X_2$).
  This set can be perfectly matched
  when $m$ is odd and can be matched except for a single face with 
  $mn/4$ particles when $m$ is even.
\qed
\end{proof}

It remains to describe the faces $\sigma$ which are in the set $P_3$.
In order to do this, we will further divide this set into three
disjoint sets, $Q_1$, $Q_2$ and $Q_3$.
\begin{itemize}
\item $Q_1$ is the set of faces $\sigma$ with signature $[3,3,\ldots,3]$.
\item $Q_2$ is the set of faces $\sigma$ for which all the intervals are
of odd length, at least one of length greater than 3, 
and such that all even positions in the intervals are non-free.
\item $Q_3$ is the set of faces $\sigma$ with signature $\emptyset$.
\end{itemize}
These sets partition the remaining faces since
all $\sigma$ with both odd and even intervals are in $P_1$ and
all $\sigma$ with only even intervals are in $P_2$.

\begin{lemma2} \label{lemma:3}
  \[
  \sum_{\sigma \in Q_1} (-1)^{|\sigma|} =
  \begin{cases}
    0 & \text{if $m \equiv_3 0$,} \\
    3 \cdot (-1)^{n/3} & \text{if $m \equiv_3 1$,} \\
    3 & \text{if $m \equiv_3 2$.} \\
  \end{cases}
  \]
\end{lemma2}

\begin{proof}
  Note that the faces with signature $[3,3,\ldots,3]$ fall into
  three equivalence classes of $\sim$.
  Let $C'_{m,n} = C_{m,n} \setminus \{ (1,3i+1) \;|\; i = 0, \ldots, n/3-1 \}$.
  We will first determine the homotopy type of $I(C'_{m,n})$.
  We then obtain the desired sum by determining
  \[
  \sum_{\sigma \in C'_{m-1,n}} (-1)^{|\sigma|+n/3}
  \]
  and multiply this by 3.

  \begin{figure}[htbp]
    \begin{center}
      \input{mtreethree.pstex_t}
    \end{center}
    \caption{Matching tree for faces $\sigma$ with $\pi(\sigma) = \{ (1, 3i+1) \;|\; i = 0, \ldots, n/3-1 \}$.}
    \label{fig:threethree}
  \end{figure}
  Figure~\ref{fig:threethree} shows a sketch of the
  matching tree construction.
  If we split on $v_1$, the subtree where $v_1$ belongs to all
  faces will be perfectly matched on $v_2$.
  We can therefore assume that $v_1$ is
  empty, and then match on $v_3$.
  In the second graph, these vertices
  have been fixed.
  We now have two vertices of degree 1, which we can match in any order
  to obtain the third graph.
  Note that after removing the first three rows and a total of
  $2n/3$ particles,
  the third graph is identical to the graph in the root.
  This provides a recursive relation.
  The base cases are when $m = 1$, $m = 2$ and $m = 3$.
  The first two of these leave a single configuration with 
  $n/3$ and $2n/3$ particles, respectively.
  For $m = 3$, the entire complex is matched.
  \[
  I(C'_{m,n}) \simeq \begin{cases}
    \bullet & \text{if $m \equiv_3 0$,} \\
    S^{2(m+1/2)n/9-1} & \text{if $m \equiv_3 1$,} \\
    S^{2(m+1)n/9-1} & \text{if $m \equiv_3 2$.} \\
  \end{cases}
  \]
  The lemma follows.
\qed
\end{proof}

\begin{conjecture} \label{conj:match}
  Let $n$ be odd.
  Then, for $\sigma \in Q_2$,
  \[
  \sum_{\tau : \pi(\tau) = \pi(\sigma)} (-1)^{|\tau|} = 0.
  \]
\end{conjecture}

We can prove the following, which shows a special case of
Conjecture~\ref{conj:match}, but also holds for even $n$.
\begin{lemma2} \label{lemma:conjpart}
  Let $\sigma \in Q_2$ and assume that for all intervals $N_i \in N(\sigma)$, it holds that $|N_i| \geq 5$.
  Then,
  \[
  \sum_{\tau : \pi(\tau) = \pi(\sigma)} (-1)^{|\tau|} = 0.
  \]
\end{lemma2}

  \begin{figure}[htbp]
    \begin{center}
      \input{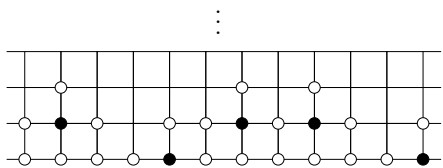}
    \end{center}
    \caption{Configurations from $Q_2$ with $n = 12$ and $\pi(\sigma) = \{5,12\}$.}
    \label{fig:confq2}
  \end{figure}

\begin{figure}[htbp]
    \begin{center}
      \input{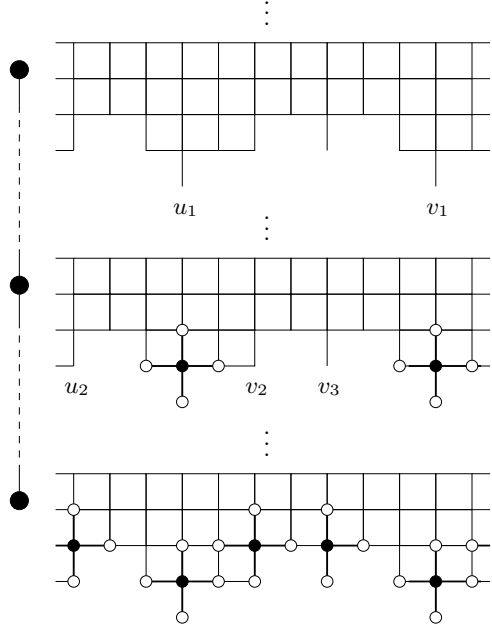}
    \end{center}
    \caption{Matching tree recursion for the configurations from Figure~\ref{fig:confq2}.}
    \label{fig:mtree6}
  \end{figure}

\begin{proof}
  Figure~\ref{fig:confq2} shows the fixed vertices in the lower part
  of a typical configuration in $Q_2$.
  In this particular example, $n = 12$ and $\pi(\sigma) = \{(1,5),(1,12)\}$.
  The same graph, when the fixed vertices are removed, is shown in
  the root of Figure~\ref{fig:mtree6}.
  We will now construct a matching tree on this graph and show that
  we can match all configurations $\tau$, with $\pi(\tau) = \pi(\sigma)$
  perfectly.

  First, we match on $u_1$ and $v_1$.
  By doing this, we obtain two new vertices of degree 1, namely $u_2$ and
  $v_2$. We also have the vertex $v_3$ which had degree 1 already in
  the root node.
  We continue by matching on $u_2, v_3$ and $v_2$ (in any order), 
  and obtain the final graph in Figure~\ref{fig:mtree6}.
  Removing all fixed vertices in this graph yields, after a rotation
  of the cylinder, the graph in the root node of height 2 less.
  This is the recursive step.
  It is easy to see that the base cases, for $m = 2$ and $m = 3$
  are perfectly matched. When $m = 1$, the set $Q_2$ is empty.

  We argue that this construction is extendable to any $\pi(\sigma)$,
  where $\sigma \in Q_2$.
  Take such a $\sigma$ and an interval $N_i \in N(\sigma)$.
  On the second row, there are $(|N_i|-3)/2$ forced particles,
  due to the fact that all even positions in $N_i$ must be non-free.
  Let $y_i$ be the last such particle, which always exists since
  $|N_i| \geq 5$.
  Then, for each $i$, $y_i$ and $x_i$ forces a degree-1 vertex.
  In the example these degree-1 vertices are $u_1$ and $v_1$.
  We start by matching on these vertices.
  In the second step, we note that for each interval, there will
  now be $(|N_i|-3)/2$ vertices of degree 1 on row 3.
  In the example these vertices are $u_2$ for the interval
  $\{(1,1),\ldots,(1,5)\}$ and $v_2$ and $v_3$ for the interval
  $\{(1,6),\ldots,(1,12)\}$.
  As in the example, when $m=2$ and $m=3$ everything can be perfectly
  matched.
  Therefore, the number of intervals does not matter, as long as
  they are all odd, and all of size greater than or equal to 5.
\qed
 \end{proof}

The final set $Q_3$ will function as an inductive step, since we have:
\begin{equation} \label{eq:q3rec}
\sum_{\sigma \in Q_3} (-1)^{|\sigma|} = \sum_{\sigma \in I(C_{m-1,n})} (-1)^{|\sigma|}.
\end{equation}

\begin{proposition2}
  Conjecture~\ref{conj:match} implies
  Conjecture~\ref{conj:cylinder}.
\end{proposition2}

\begin{proof}
  Let $n$ be odd.
  We know from Lemma~\ref{lemma:perfect} that the alternating sum over
  $P_1$ is 0.
  Since $m$ is odd, we further know that $P_2$ is empty.
  For $P_3$, we use Lemma~\ref{lemma:3} for $Q_1$ and 
  Conjecture~\ref{conj:match} for $Q_2$.
  Finally, as noted in $(\ref{eq:q3rec})$, we get an
  inductive step from $Q_3$.
  \begin{multline*}
    Z(C_{m,n}) = \sum_{\sigma \in P_1} (-1)^{|\sigma|} + \sum_{\sigma \in P_2} (-1)^{|\sigma|} + \sum_{\sigma \in P_3} (-1)^{|\sigma|} = \\
    = 0 + 0 + 
    \sum_{\sigma \in Q_1} (-1)^{|\sigma|} + 
    \sum_{\sigma \in Q_2} (-1)^{|\sigma|} +
    \sum_{\sigma \in Q_3} (-1)^{|\sigma|} = \\
    = Z(C_{m-1,n}) +
    \begin{cases}
      0 & \text{if $m \equiv_3 0$ and $3|n$,} \\
      -3 & \text{if $m \equiv_3 1$ and $3|n$,} \\
      3 & \text{if $m \equiv_3 2$ and $3|n$.} \\
    \end{cases}
  \end{multline*}
  The base case is $Z(C_{1,n})$ which is the alternating sum
  of the independence complex on an $n$-cycle.
  Its values are given by the following expression.
  \[
  Z(C_{1,n}) = \begin{cases}
    -2 & \text{if $n \equiv_3 0$,} \\
    (-1)^{(n-1)/3} & \text{if $n \equiv_3 1$,} \\
    (-1)^{(n+1)/3} & \text{if $n \equiv_3 2$.} \\
  \end{cases}
  \]
  Since $n$ is odd, this can be simplified to $Z(C_{1,n}) = -2$ if $3|n$ and $Z(C_{1,n}) = 1$ otherwise.
  Together with the inductive step, this shows the proposition.
\qed
\end{proof}

When $n = 3$, the set $Q_2$ is empty and for $n = 5, 7$ and 9, 
Lemma~\ref{lemma:conjpart} covers all the cases of 
Conjecture~\ref{conj:match}.
Therefore, Conjecture~\ref{conj:cylinder} holds for these values of $n$.

\begin{figure}[htbp]
  \begin{center}
    \input{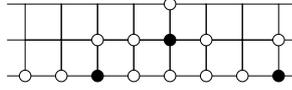}
  \end{center}
  \caption{Counter-example to Conjecture~\ref{conj:match} when $n$ is even.}
  \label{fig:counter}
\end{figure}

\medskip
Conjecture~\ref{conj:match} does not hold for all even $n$.
A counterexample is given in Figure~\ref{fig:counter}, where
$m = 3, n = 8$ and $\pi(\sigma) = \{(1,3),(1,8)\}$.
A straightforward application of the matching tree construction
produces a matching on the faces $\tau$ with 
$\pi(\tau) = \pi(\sigma)$ which leaves a single face with six
particles. Thus, 
$\sum_{\tau : \pi(\tau) = \{3,8\}} (-1)^{|\tau|} = 1$ in this case.

Still, we can determine $Z(C_{m,n})$ for some small, even $n$.
When $n = 6$, both $P_1$ and $Q_2$ are empty.
Lemma~\ref{lemma:p2}, Lemma~\ref{lemma:3} and 
$(\ref{eq:q3rec})$ gives
\[
\sum_{\sigma \in I(C_{6,n})} (-1)^{|\sigma|} =
- \sum_{\sigma \in I(C_{m-1,n})} (-1)^{|\sigma|} +
\begin{cases}
  0+2 \cdot (-1)^{m/2} & m \equiv_6 0 \\
  3+0 & m \equiv_6 1 \\
  3+2 \cdot (-1)^{m/2} & m \equiv_6 2 \\
  0+0 & m \equiv_6 3 \\
  3+2 \cdot (-1)^{m/2} & m \equiv_6 4 \\
  3+0 & m \equiv_6 5 \\
\end{cases}
\]
Therefore, $Z(C_{m,6})$ is the following periodic sequence,
starting from $m=1$.
The period is 12.
\[
{\bf 2}, -1, 1, 4, -1, -1, 4, 1, -1, 2, 1, 1, {\bf 2}, -1, 1, \ldots
\]


\begin{table}[htbp]
  \begin{center}
    \begin{tabular}{|c|r||*{11}{r|}}
      \hline
      \multicolumn{2}{|c||}{\multirow{2}{*}{$\displaystyle\sum_{\sigma \in Q_2} (-1)^{|\sigma|}$}} & \multicolumn{11}{|c|} m \\
      \cline{3-13}
      \multicolumn{2}{|c||}{}&\makebox[14pt][r]{1}&\makebox[14pt][r]{2}&\makebox[14pt][r]{3}&\makebox[14pt][r]{4}&\makebox[14pt][r]{5}&\makebox[14pt][r]{6}&\makebox[14pt][r]{7}&\makebox[14pt][r]{8}&\makebox[14pt][r]{9}&\makebox[14pt][r]{10}&\makebox[14pt][r]{11} \\
      \hline
      \hline
      \multirow{11}{*}{\,\,\, n \,\,\,}
        & 2 &   &   &    &    &   &   &    &    &   &   &    \\
      \cline{2-13}
        & 4 &   &   &   &    &   &    &   &    &   &   &   \\
      \cline{2-13}
        & 6 &    &   &    &    &   &   &    &    &   &    &    \\
      \cline{2-13}
        & 8 &   &    &  8 &    &  8 &    &  8 &  -8 & 8 &  -8 & 8 \\
      \cline{2-13}
        &10&   &   &    &    &  10 &   &    &    &-10 &   &    \\
      \cline{2-13}
        &12&    &    &   &    &   &    &  12 &    &   &   &   \\
      \cline{2-13}
        &14&   &   &    &    &   & 14 &    &    & 14 &   & 14 \\
      \cline{2-13}
        &16&   &    &  8 &    &  8 &    &  8 & 24 & 8 &  -8 & 24 \\
      \cline{2-13}
        &18&    &   &    &    &   &   & 18 &    &   & 36 &    \\
      \cline{2-13}
        &20&   &    &   &    &  10 &    &   &    & 50 &   &   \\
      \cline{2-13}
        &22&   &   &    &    &   &   &    & 22 &   &   & 88 \\
      \cline{2-13}
        &24&    &    &  8 &    &  8 &    & 20 &  -8 & 8 & 64 & 8 \\
      \hline
    \end{tabular}
  \end{center}
  \caption{Alternating sums over $Q_2$ for small $m$ and $n$.}
  \label{tab:sqcylminus}
\end{table}

In Table~\ref{tab:sqcylminus} we have taken the values of $Z(C_{m,n})$
from Table~\ref{tab:sqcyltrans} and removed the part
which can be explained using
Lemma~\ref{lemma:p2}, Lemma~\ref{lemma:3} and $(\ref{eq:q3rec})$
similarly to what we did in the previous example.
What remains in the table is therefore
the alternating sums over the faces in
$Q_2$, or alternatively over those faces in
 $Q_2$ which are not matched by Lemma~\ref{lemma:conjpart}.
These are still unexplained, but
several interesting patterns seem to emerge.

\section{Transfer matrices} \label{sec:transfermatrices}

In \cite{BMLN07}, Bousquet-Mélou, Linusson and Nevo determined the partition function of the hard particle model on a parallelogram shaped region $P_{m,n}$ of the square grid. We describe it by letting $(a,b) \in V(P_{m,n})$ if
\[
1 \leq a \leq m \quad \text{ and } \quad 1-a < b \leq 1-a + n.
\]
\begin{figure}[htbp]
  \begin{center}
    \input{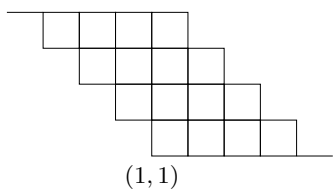}
  \end{center}
  \caption{The graph $P_{5,6}$.}
  \label{fig:parallel}
\end{figure}

Figure~\ref{fig:parallel} shows the region $P_{5,6}$.
The skewed border allowed them to find a Morse-matching 
for the complex $I(P_{m,n})$ and determine its homotopy type.
They further realised that this could be used to produce 
parts of the spectrum of the transfer matrix for the cylinder $C_{m,n}$.
In this section, we will review this idea and show how we can
complete the case for $m = 4$.
We also find the missing eigenvalue when $m = 6$.
Unfortunately, we were not able to generalise our findings
to new, unknown spectra.

\medskip
We start by defining the two transfer matrices of interest.
For a fixed $m$, the matrix $T(m)$ describes the configurations
on the region $P_{m,n}$ and
$T'(m)$ describes the configurations on the region $S_{m,n}$.
Throughout this section, $i$ will denote the square root of $-1$.

Let $T(m)$ be the $2^m \times 2^m$-matrix, indexed by
subsets of $[m]$ on both rows and columns, and defined as follows.
\begin{equation} \label{eq:tp}
  T(m)_{A,B} = \begin{cases}
    i^{|A|+|B|} & \text{if } A \cap B = A \cap B' = \emptyset, \\
    0 & \text{otherwise,}
  \end{cases}
\end{equation}
where $B' = \{ j + 1 \;|\; j \in B \cap [m-1] \}$.

Let $T'(m)$ be the $F_{m+1} \times F_{m+1}$-matrix, indexed
by independent sets on a path with $m$ vertices, and defined as follows.
\begin{equation} \label{eq:tpp}
  T'(m)_{A,B} = \begin{cases}
    i^{|A|+|B|} & \text{if } A \cap B = \emptyset, \\
    0 & \text{otherwise.}
  \end{cases}
\end{equation}
We here follow \cite{BMLN07} and
in contrast to the example in Section~\ref{subsec:hardp}
distribute the weights of the entries in row $A$, column $B$
and row $B$, column $A$ evenly.

As was illustrated in Section~\ref{subsec:hardp},
we can determine the values of $Z(C_{m,n})$ by taking
the trace of the $n$th power of either one of these matrices.
Therefore, the following relation between the
generating functions of $Z(C_{m,n})$ holds.
\begin{equation} \label{eq:spectra}
  \textrm{tr}\left(1 - t T(m)\right)^{-1} - 2^m = 
  \sum_{n \geq 1} Z(P_{m,n}) t^n =
  \textrm{tr}\left(1 - t T'(m)\right)^{-1} - F_{m+1}
\end{equation}
We conclude that the spectra of $T(m)$ and $T'(m)$ are equal,
apart from the multiplicity of the null eigenvalue
due to the difference in size.

Let $A, B \subseteq [m]$ and assume that we prescribe the particles
given by $A$ and $B$ on the left and right boundary of $P_{m,n}$.
That is, when $j \in A$, then the $j$th particle from the top on the
left boundary is present, and otherwise it is missing in each configuration
and similarly for $B$ on the right boundary.
This induces a subgraph of $P_{m,n}$ which we will denote by $P_{m,n}(A,B)$.

Let $G_{A,B}(t)$ denote the element is row $A$, column $B$ of 
$(1-tT(m))^{-1}$.
Then, $G_{A,B}(t)$ is the generating function for configurations
on $P_{m,n}(A,B)$, weighted as given by $(\ref{eq:tp})$
and has the following expression.
\begin{equation} \label{gdef}
  G_{A,B}(t) = \delta_{A,B} + (-i)^{|A|+|B|} \sum_{n\geq 1} Z(P_{m,{n+1}}(A,B)) t^n
\end{equation}

We can then split the left-hand side of $(\ref{eq:spectra})$ into
\begin{equation} 
  \textrm{tr}\left(1- t T(m)\right)^{-1} = \sum_{A \subseteq [m]} G_{A,A}(t).
\end{equation}

The observation made in \cite{BMLN07} was that parts of the spectra
can be recovered by finding $G_{A,A}(t)$ for some $A$ and
that this can be done in the model for $P_{m,n}$ which is more
easily determined than that of $C_{m,n}$.

The characteristic polynomial of $T'(4)$ is determined both in
\cite{BMLN07} and \cite{Jonsson06} and
is given by $(1-t+t^2) (1-t^2) (1-t^4)$.
Let $G_{A,B}^{(N)}(t)$ denote the first $N$ terms in the expansion of
$(\ref{gdef})$.
By finding a Morse-matching for the region in Figure~\ref{fig:parallel}, 
Bousquet-Mélou, Linusson and Nevo effectively calculated 
$G_{\emptyset,\emptyset}(t)$.
When $m = 4$, they found that

\begin{equation} \label{eq:emptyempty4}
  G_{\emptyset,\emptyset}(t) = \frac{1}{1-t+t^2}.
\end{equation}
This determined

This recovers one of the factors of $T'(4)$.
Using this case, and a recursion for $Z\left(P_{m,n}(\{2,3\},\{2,3\})\right)$,
they could then determine the function $G_{\{2,3\},\{2,3\}}(t)$.
\begin{equation}
  G_{\{2,3\},\{2,3\}}(t) = \frac{1+t^2+t^3}{(1+t^3)(1-t^4)}
\end{equation}

We will now determine the remaining cases for $m=4$.
We start with the following lemma.
\begin{lemma2} \label{lem:match}
  Assume that $j-1, j, j+1 \in A$ for some $j$ and let $B = A \setminus \{j\}$.
  Then,
  \[
  G_{A,A}(t) + G_{B,B}(t) = 2.
  \]
\end{lemma2}

\begin{figure}[htbp] 
  \begin{center}
    \input{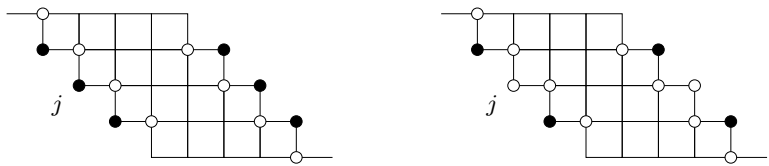}
  \end{center}
  \caption{Two regions differing only in the position $j$.}
  \label{fig:equal}
\end{figure}

\begin{proof}
  Figure~\ref{fig:equal} illustrates that the regions remaining after fixing
  the boundaries to $A$ and $B$ respectively are equal.
  Since the addition of two elements does not change the parity, 
  we have $Z\left(P_{m,n}(A,A)\right) = Z\left(P_{m,n}(B,B)\right)$.
  Since furthermore $(-i)^{|A|+|A|} = - (-i)^{|B|+|B|}$,
  all terms except the constants cancel and
  we have $G_{A,A}(t) + G_{B,B}(t) = \delta_{A,A} + \delta_{B,B} = 2$.
\qed
\end{proof}

Using Lemma~\ref{lem:match} we can avoid calculating certain sets by
choosing a partial matching on $2^{[4]}$ which satisfies the condition 
in the lemma.
We will in this way match $\{1,2,3\}$ with $\{1,3\}$, $\{2,3,4\}$ 
with $\{2,4\}$ and $\{1,2,3,4\}$ with $\{1,3,4\}$.
This leaves only $\{1,2,4\}$ unmatched among the potential candidates
for Lemma~\ref{lem:match}.
We can produce its generating function by noting that 
$G_{[4],[4]}(t) = 1 + t^2 G_{\emptyset,\emptyset}(t)$.
Lemma~\ref{lem:match} applied to $A = [4]$ and $i = 3$ 
gives us
\begin{equation}
  G_{\{1,2,4\},\{1,2,4\}}(t) = 1 - t^2 G_{\emptyset,\emptyset}(t) = 
  \frac{1-t}{1-t+t^2}.
\end{equation}

The sets that remain are the singleton sets, $\{1,2\}$, $\{3,4\}$ and $\{1,4\}$.
By symmetry, we can restrict our calculations to the sets 
$\{3\}$, $\{4\}$, $\{1,2\}$ and $\{1,4\}$.
We present the calculations for $\{4\}$ and $\{1,4\}$.
\medskip

We start with the case $A = \emptyset$ and $B = \beta$, 
where $\beta$ is an arbitrary subset of $[m]$.
When $j \geq 4$, we have the recursive relation
$Z\left(P_{4,j}(\emptyset,\beta)\right) = -Z\left(P_{4,j-3}(\emptyset,\beta)\right)$.
Thus, $G_{\emptyset,\beta}(t) = G^{(4)}_{\emptyset,\beta}(t) - t^3 G_{\emptyset,\beta}(t)$ which gives us the following general expression.
\begin{equation} \label{eq:0beta}
  G_{\emptyset,\beta}(t) = \frac{G^{(4)}_{\emptyset,\beta}(t)}{1+t^3}
\end{equation}

The recursive relation which was used to produce $(\ref{eq:0beta})$, 
and those that follow, are obtained using simple
matching trees.
We present some, but not all trees here.
\begin{figure}[htbp]
  \begin{center}
    \input{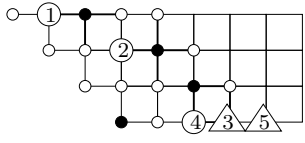}
  \end{center}
  \caption{The large circles indicate match nodes and the
    large triangle indicate split nodes.
    The numbers give the order in which the nodes are matched or split.
  At 5, the two branches of the recursion are obtained.}
  \label{fig:rec1}
\end{figure}
Figure~\ref{fig:rec1} illustrates a matching tree when $A = \{4\}$.
This gives us the following recursive relation for the case $A = B = \{4\}$.
\begin{equation} \label{eq:rec44}
  Z\left(P_{4,j}(\{4\},\{4\})\right) = Z\left(P_{4,j-4}(\{4\},\{4\})\right) + Z\left(P_{4,j-4}(\emptyset,\{4\})\right)
\end{equation}
For the second term, we calculate $G^{(4)}_{\emptyset,\{4\}}(t) = -it$, 
which, using $(\ref{eq:0beta})$, gives
$G_{\emptyset,\{4\}}(t) = - it / (1+t^3)$.
Plugging this term into $(\ref{eq:rec44})$ gives the expression
\begin{equation}
  G_{\{4\},\{4\}}(t) = G_{\{4\},\{4\}}^{(5)}(t) + t^4 \left(G_{\{4\},\{4\}}(t)-1\right) + i t^4 G_{\emptyset, \{4\}}(t),
\end{equation}
which, after calculating $G_{\{4\},\{4\}}^{(5)}(t) = 1 + t^4$,
gives the final function for $A = B = \{4\}$, and by symmetry for
$A = B = \{1\}$.
\begin{equation} \label{eq:44}
G_{\{4\},\{4\}}(t) = G_{\{1\},\{1\}}(t) = \frac{1+t^3+t^5}{(1+t^3)(1-t^4)}
\end{equation}

For $A = B = \{1,4\}$, we have the recursion 
\begin{equation} \label{eq:rec1414}
  Z(P_{4,j}(A,A)) = Z(P_{4,j-1}(\emptyset,A)) + Z(P_{4,j-4}(A,A).
\end{equation}
We proceed by summing over $j$ and calculating the base cases.
\begin{multline}
G_{\{1,4\},\{1,4\}}(t) \\ = G_{\{1,4\},\{1,4\}}^{(5)}(t) - t \left( G_{\emptyset,\{1,4\}}(t) - G_{\emptyset,\{1,4\}}^{(4)}(t) \right) + t^4 \left( G_{\{1,4\},\{1,4\}}(t)-1 \right) \\
= \left( 1 + t^2 + t^4 + \frac{t^2}{1+t^3} - t^2 - t^4 \right) / (1-t^4)
\end{multline}
In the last equality, we have calculated the polynomials
\[
G_{\{1,4\},\{1,4\}}^{(5)}(t) = 1 + t^2 + t^4 \text{ and }
G_{\emptyset,\{1,4\}}^{(4)}(t) = -t,
\]
and used the expression for $G_{\emptyset, \{1,4\}}(t)$ 
in $(\ref{eq:0beta})$.
It turns out that $G_{\{1,4\},\{1,4\}}(t) = G_{\{2,3\},\{2,3\}}(t)$.

\smallskip
The final two cases, when $A = B = \{3\}$ and $A = B = \{1,2\}$,
are determined similarly, and we present them without any
explicit calculations.
\begin{equation} \label{eq:33}
  G_{\{2\},\{2\}}(t) = G_{\{3\},\{3\}}(t) = \frac{1+t^3+t^5}{(1+t^3)(1-t^4)}
\end{equation}
Note that $(\ref{eq:33})$ implies that all the singleton sets 
have the same generating function.
\begin{equation} \label{eq:1212}
G_{\{1,2\},\{1,2\}}(t) = G_{\{3,4\},\{3,4\}}(t) = 
1 - \frac{t^5}{(1+t^3)(1-t^4)}.
\end{equation}

We note that while all factors of the characteristic polynomial
show up in the individual terms, we can not deduce their
multiplicities until we have taken the sum over all sets.
Expanding this sum into partial fractions, 
we arrive at the following expression.

\begin{multline*}
  \sum_{A \subseteq [4]} G_{A,A}(t) =
  G_{\emptyset,\emptyset}(t) + G_{\{1,2,4\},\{1,2,4\}}(t) + 6 + \\
  + 4 G_{\{4\},\{4\}}(t) + 2 G_{\{2,3\},\{2,3\}}(t) + 2 G_{\{1,2\},\{1,2\}}(t)
  \\
  = 8 + \frac{1}{1-it} + \frac{1}{1+it} + \frac{2}{1+t} + \frac{2}{1-t} + \frac{1}{1-\alpha t} + \frac{1}{1-{\bar \alpha} t},
\end{multline*}
where $\alpha = (1+i\sqrt{3})/2$ and ${\bar \alpha}$ are the two roots of $1-t+t^2$.
As expected, we recover the roots and multiplicities of the
characteristic polynomial of $T(4)$, which apart from the 8 zeroes
equals the characteristic polynomial of $T'(4)$.

\smallskip
Next, we determine the factor of the characteristic polynomial
for $m = 6$ which does not follow from the case $A = B = \emptyset$.
The characteristic polynomial for $m=6$ is given by
$(1-t^4)^2 (1-t^{14})/(1+t)$ \cite{BMLN07,Jonsson06}
and the function $G_{\emptyset,\emptyset}(t)$ follows from
the general result in \cite{BMLN07}.
\begin{equation}
  G_{\emptyset, \emptyset}(t) = \frac{1}{1-t^{14}}\left(1+t+t^3 \frac{1-t^{12}}{1-t^3} \right)
\end{equation}
This leaves the factor $1-t^4$ unexplained.
We show that this factor appears in the generating function 
$G_{A,A}(t)$ for $A = \{2,3,6\}$.

\smallskip
As in the previous case, we start by deriving an expression for
$G_{\emptyset,\beta}(t)$ for an arbitrary $\beta$.
The matching leading up to the recursion contains more steps,
but is straightforward using matching trees.
The final formula is
\begin{equation}
  G_{\emptyset,\beta}(t) = \frac{G_{\emptyset,\beta}^{(15)}(t)}{1-t^{14}}.
\end{equation}

We will also make use of the recursive relation
\begin{equation}
  Z\left(P_{6,j}(\{2\},\{2,3,6\})\right) = 
  -Z\left(P_{6,j-6}(\emptyset,\{2,3,6\})\right),
\end{equation}
from which we derive that
\begin{multline}
  G_{\{2\},\{2,3,6\}}(t) = G_{\{2\},\{2,3,6\}}^{(7)}(t) +i t^6 G_{\emptyset,\{2,3,6\}}(t) \\
  = t^4 +i t^6 \left( \frac{-i(t+t^4+t^{12})}{1-t^{14}} \right)
  = \frac{t^4 \left(1 + t^3 + t^6 \right)}{1-t^{14}}.
\end{multline}

\begin{figure}[htbp]
  \begin{center}
    \input{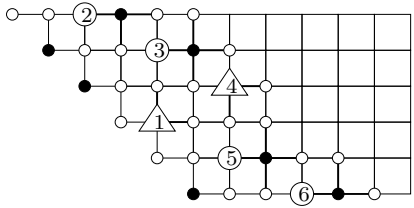}
  \end{center}
  \caption{A matching for $P_{6,j}(\{2,3,6\},\{2,3,6\})$.}
  \label{fig:rec2}
\end{figure}

In the last equality we have used that
$G_{\emptyset,\{2,3,6\}}^{(15)}(t) = -i \left( t + t^4 + t^{12} \right)$.
The main recursion for $A = \{2,3,6\}$ is given by the following
relation.
\begin{multline} \label{eq:6}
  Z\left(P_{6,j}(A,A)\right) = 
  Z\left(P_{6,j-4}(A,A)\right)
  - Z\left(P_{6,j-1}(\emptyset,A)\right)
  - Z\left(P_{6,j-4}(\{2\},A)\right)
\end{multline}
The matching tree needed to produce this relation is shown in
Figure~\ref{fig:rec2}.
The terms $Z\left(P_{6,j-1}(\emptyset,A)\right)$ and
$Z\left(P_{6,j-4}(\{2\},A)\right)$ in the right-hand side of 
$(\ref{eq:6})$ are obtained in nodes 1 and 4, respectively.
After summation over $i$ we have
\begin{multline} \label{eq:g236}
  G_{\{2,3,6\},\{2,3,6\}}(t) = G_{\{2,3,6\},\{2,3,6\}}^{(5)}(t) + t^4 \left( G_{\{2,3,6\},\{2,3,6\}}(t) - 1 \right) \\
  -it \left( G_{\emptyset,\{2,3,6\}}(t) - G_{\emptyset,\{2,3,6\}}^{(4)}(t) \right) + t^4 G_{\{2\},\{2,3,6\}}(t)
\end{multline}

Calculations of initial terms in the relevant series
provides the polynomials
\[
G_{\{2,3,6\},\{2,3,6\}}^{(5)}(t) = 1-t^2+t^4 \text{ and }
G_{\emptyset,\{2,3,6\}}^{(4)}(t) = -i t,
\]
so $(\ref{eq:g236})$ evaluates to
\begin{equation}
  G_{\{2,3,6\},\{2,3,6\}}(t) = 
  \frac{(1+t-t^5-t^6-t^7+t^{11}+t^{12})(1-t)}{(1-t^{14})(1-t^4)}.
\end{equation}
This recovers the missing factor $1-t^4$.


\bigskip
{\bf Acknowledgements }
\noindent
The author wishes to thank his supervisor Svante Linusson, who
provided many ideas and suggestions to the work presented in this paper
and who also
patiently read and commented on countless early manuscripts.


\bibliographystyle{plain}

\end{document}